\documentclass[11pt]{article}
\usepackage{jheppub,amsmath,  amssymb,slashed,url,bm,textgreek,upgreek}
\usepackage{graphicx}
\usepackage{epstopdf}
\newtheorem{theorem}{Theorem}[section]
\newtheorem{prop}[theorem]{Proposition}
\newtheorem{lemma}[theorem]{Lemma}

\newtheorem{definition}[theorem]{Definition}
 
\def\KW{\bf{KW}}

\def\epsilon{\varepsilon}
\def\tt{{\frak t}}

\def\Tr{{\mathrm{Tr}}}

\def\X{{\mathcal X}}

\def\s{{s}}

\def\p{{\eurm p}}

\def\be{\begin{equation}}
\def\ee{\end{equation}}

\def\X{{\eusm X}}

\def\hat{\widehat}

\def\tilde{\widetilde}

\def\frak{\mathfrak}

\def\D{{\mathcal D}}

\def\V{{\mathcal V}}

\def\Bbb{\mathbb}

\def\d{{\mathrm d}}

\def\R{{\mathbb R}}
\def\C{{\mathbb C}}

\def\D{{\mathcal D}}
\def\[{\bigl [}

\def\]{\bigr ]}

\def\Z{{\mathbb Z}}

\def\ad{{\mathrm{ad}}}
\def\r{{\mathrm{tr}}}

\def\V{{\mathcal V}}

\def\r{\rangle}

\def\tilde{\widetilde}

\def\bar{\overline}

\def\Y{{\eusm Y}}

\font\teneurm=eurm10 \font\seveneurm=eurm7  \font\fiveeurm=eurm5
\newfam\eurmfam
\textfont\eurmfam=\teneurm \scriptfont\eurmfam=\seveneurm
\scriptscriptfont\eurmfam=\fiveeurm
\def\eurm#1{{\fam\eurmfam\relax#1}}
\font\teneusm=eusm10 \font\seveneusm=eusm7 \font\fiveeusm=eusm5
\newfam\eusmfam
\textfont\eusmfam=\teneusm \scriptfont\eusmfam=\seveneusm
\scriptscriptfont\eusmfam=\fiveeusm

\font\tencmmib=cmmib10 \skewchar\tencmmib='177
\font\sevencmmib=cmmib7 \skewchar\sevencmmib='177
\font\fivecmmib=cmmib5 \skewchar\fivecmmib='177
\newfam\cmmibfam
\textfont\cmmibfam=\tencmmib \scriptfont\cmmibfam=\sevencmmib
\scriptscriptfont\cmmibfam=\fivecmmib

\def\del{\partial}
\def\calA{\mathcal A}
\def\calO{\mathcal O}
\def\calL{\mathcal L}
\def\LKW{\mathcal L}
\newcommand{\wfa}{\mathrm{wf}}
\newcommand{\ff}{\mathrm{ff}}
\newcommand{\iie}{\mathrm{iie}}
\newcommand{\diag}{\mathrm{diag}}
\newcommand{\diagiie}{\mathrm{diag}_{\iie}}
\newcommand{\fdiag}{\mathrm{fdiag}}
\newcommand{\ie}{\mathrm{ie}}
\newcommand{\calJ}{\mathcal J}
\newcommand{\calE}{\mathcal E}
\newcommand{\calU}{\mathcal U}
\newcommand{\calV}{\mathcal V}
\newcommand{\calD}{\mathcal D}
\newcommand{\calC}{\mathcal C}
\newcommand{\RR}{\mathbb R}
\newcommand{\wh}{\widehat}

\def\R{{\Bbb R}}
\def\d{{\mathrm d}}
\def\D{{\mathcal D}}
\def\C{{\Bbb C}}
\def\hD{\widehat\Delta}
\def\uphi{\text{\textphi}}

\def\X{{\mathcal X}}
\def\Y{{\mathcal Y}}
\def\Z{{\mathcal Z}}
\def\Tr{{\mathrm{Tr}}}

\def\V{{\mathcal V}}
\def\r{{\frak r}}
\def\s{{\frak s}}
\def\bar{\overline}

\title{The KW Equations and the Nahm Pole Boundary Condition with Knots} 

\author{Rafe Mazzeo$^a$}
\affiliation{$^a$Department of Mathematics, Stanford University, Stanford, CA 94305} 

\author{and Edward Witten$^{b}$}
\affiliation{$^{b}$School of Natural Sciences, Institute for Advanced Study,\\ 1 Einstein Drive, Princeton, NJ 08540 USA}
\abstract{It is conjectured that the coefficients of the Jones polynomial can be computed by counting solutions of the KW equations on a four-dimensional half-space,
with certain boundary conditions that depend on a knot.  The boundary conditions are defined by a ``Nahm pole'' away from the knot with a further singularity along
the knot.  In a previous paper, we gave a precise formulation of the Nahm pole boundary condition in the absence of knots; in the present paper, we do this in the more
general case with knots included.  We show that the KW equations with generalized Nahm pole boundary conditions are elliptic, and that the solutions are polyhomogeneous
near the boundary and near the knot, with exponents determined by solutions of appropriate indicial equations.  This involves the analysis of a ``depth two incomplete iterated edge operator.''    As in
our previous paper, a key ingredient in the analysis is a convenient new Weitzenb\"ock  formula that is well-adapted to the specific problem.}

\begin{document} \maketitle

\section{Introduction}

\def\emptyset{\slashed{\mathrm{O}}}

In our previous paper \cite{MW}, we investigated some of the mathematical underpinnings of the Kapustin-Witten (KW) equations on a compact oriented 
Riemannian four-manifold $(X,g)$ with boundary with Nahm pole boundary condition at $\del X = W$. These are equations for a pair $(A, \phi)$ where $A$ is a 
connection on a $G$-bundle $E\to M$ and $\phi$ is a 1-form on $X$ valued in the adjoint bundle $\ad(E)$, and take the form\footnote{A more general version of the
equations depends on a real parameter $t$.  See for example eqn. (2.11) in \cite{MW}.  We here set $t=1$ for simplicity.   The results in this paper are expected to have analogs for generic
$t$, but at any rate, being true at $t=1$, they certainly hold in an open set containing $t=1$.}
\begin{equation}
\label{KW} 
\begin{aligned}
F_A -\phi\wedge\phi+\star\,\d_A\phi &=0 \\
\d_A\star \phi & = 0,
\end{aligned}
\end{equation}
where $\star$ is the Hodge star and $\d_A=\d+[A,\cdot]$ is the extension of the connection to a map between differential forms 
of any degree. Informally, the boundary condition at $W$ states that $A$ is continuous up to $W$ and that
\begin{equation}\label{npbc} 
\phi \sim \frac{1}{y} \phi_\varrho;
\end{equation}
here $y$ is a boundary defining function, $\varrho$ is a principal representation $\frak{su}(2) \to \frak{g}$ and $\phi_\rho$ is an associated injective 
bundle map $TW \to \ad(E)$. 

In this paper we generalize this by considering a situation where the boundary $W$ also carries an embedded one-manifold $K$, 
i.e., a knot or link, along which there is a more subtle singularity.  This singularity depends on the choice, for each component in $K$, of a dominant weight (or 
equivalently, an irreducible representation) of the Langlands or GNO dual group $G^\vee$ to $G$.   The motivation for this generalization is that it is conjectured
\cite{W5BK} that the coefficients of the Jones polynomial of a knot can be computed by counting solutions of the KW equations on a half-space in $\R^4$  with the generalized Nahm pole boundary conditions
that we will explore.
The Jones polynomial is a Laurent series $J(q)=\sum_n a_n q^n$, and the conjecture is that $a_n$ is an algebraic count of the number of solutions of the KW equations
with second Chern class $n$.  For a more precise statement of the conjecture, see \cite{W5BK}.

Our goal here is to find analogues of the various results 
in \cite{MW} in the presence of this extra structure.  In contrast to the previous paper (where arbitrary embeddings were considered), we assume that the Nahm pole singularity at a generic point of $W$ 
is associated to  a regular embedding $\varrho:\frak{su}(2)\to\frak g$.
The generalization to arbitrary $\varrho$ presents technical difficulties of that can be handled by known methods in the absence of knots, but the corresponding extension
when $K \neq \emptyset$, while fully expected to hold, will perhaps require significant extra work.  The key difference is that when $\varrho$ is regular,
the model solution around the knot is unique, and this has ramifications throughout both the linear and nonlinear analysis.\footnote{In the nonregular case, there is a family of model solutions depending on some parameters.  This family (which depends on $\varrho$ and 
also on a weight of the dual group) has not yet been studied carefully, but this would certainly be necessary to generalize the analysis here to that case. 
 For $G$ of rank bigger than 1
but $\varrho$ still regular, the model solution is still unique (see \cite{Mikhaylov} for a construction of these solutions)
so an extension of our considerations  in that direction would be potentially more straightforward.   
If $\varrho$ is nonregular, the application of the Nahm pole boundary condition to knot theory is not fully understood, but in some cases
leads to quantum knot invariants associated to a supergroup that has $G$ as a subgroup   \cite{MiW}.
For example, for $G=\mathrm{SU}(2)$ and $\varrho=0$, the supergroup is $\mathrm{SU}(2|1)$.}  In addition to taking $\varrho$ to be regular, we will in this paper choose
$G=\mathrm{SU}(2)$ or $\mathrm{SO}(3)$,  which means that $\varrho$ is either identically zero (meaning that there is no singularity at a generic boundary
point, although there is still a singularity along $K$) or else is regular.   The regular case, which we consider here, is the case that 
 $\varrho:\frak{su}(2)\to \frak{g}$ is an isomorphism.

The starting point for this generalization is a description of the new boundary condition along $K$.  The main ingredient in this is
a model solution $(A^m, \phi^m)$ , a derivation of which appears below, on the model half-space $\RR^4_+$ with singularity
along a straight line $\RR \subset \RR^3 = \del \RR^4_+$.  This solution is translation invariant, and hence solves a reduced equation on $\RR^3_+$; 
it has a Nahm singularity along $\del \RR^3_+ \setminus \{0\}$ and in addition blows up like the inverse of the polar distance to the origin. 
This pair can be transferred to each fiber of the inward-pointing normal bundle to $K$ in $M$, $N^+K$, resulting in an approximate solution $(A^K, \phi^K)$ 
defined in a neighborhood of $K$. We then seek solutions of \eqref{KW} with a Nahm pole singularity as described above along $W \setminus K$ and with
$$
(A, \phi) \sim (A^K, \phi^K)
$$
near $K$.  One of the key tasks here is to formulate this precisely with precise rates of decay for the difference between the two sides.

One part of the story in \cite{MW} revolves around a new Weitzenb\"ock  formula specially adapted to this singular boundary condition. 
Another important step there is the determination of the formal rate of decay of homogeneous solutions of the linearization of the KW equations,
which is a lengthy and essentially algebraic calculation. This is accompanied by a regularity theorem which shows that these formal
decay rates correspond to terms in a polyhomogeneous expansion of suitably gauged solutions to \eqref{KW} at the boundary. This is necessary 
for the precise formulation of the Nahm pole boundary condition, as well as for justifying the various manipulations and calculations. 
The other main result in \cite{MW} is the calculation of the index of this operator.

In the present paper we obtain analogues of each of these results in the presence of a knot in the boundary. There are some important new features.  
First, determination of the formal rates of decay of solutions of the linearization is no longer a purely algebraic problem. These rates, also known as 
indicial roots, are calculated in this case in terms of the eigenvalues of an induced elliptic operator on the unit half-sphere in each fiber of $N^+K$.
Thus we can no longer give their precise values, but must at least produce lower bounds for the indicial roots near $0$.  
Another key difference is that in this case the linearized gauged KW operator now has a more complicated singular structure (technically it is 
a ``depth-two incomplete iterated edge operator''). The analysis needed to understand its Fredholmness and regularity properties of its solutions 
must take into account the Nahm pole and knot singularities separately. This is done using tools from geometric microlocal analysis. 
These constructions are generalizations of ones in the `pseudodifferential edge calculus' \cite{Ma}, as employed in \cite{MW}, but there 
are a number of new features here, including a structured iteration to obtain the sharp regularity theorem.

\section{The Model Solution at a Knot}\label{gennps} 
The first task is to give an explicit description for the leading order singularity imposed on solutions near $K$.  We recall an explicit formula for this in the model
case where $X$ is a Euclidean half-space $\RR^4_+$ and $K$ is a straight line in the boundary, $\RR \subset \RR^3  = \del \RR^4_+$.  The calculation below is
taken from \cite{W5BK}, starting in Section 3.6.2, and is based on an alternate expression for the KW equations.  As explained in the introduction, we assume that
$G = \mathrm{SU}(2)$ or $\mathrm{SO(3)}$.  
We use linear coordinates $\vec{x} = (x^1, x^2, x^3)$ along the boundary of $\RR^4_+$ and a normal coordinate  $x^4 = y \geq 0$, and we fix the Euclidean metric $g = \sum_{j=1}^4 (\d x^j)^2$. 

We seek a model solution to the KW equations on this half-space which is invariant under translations in $x^1$ and with the property that $A_1=\phi_4=0$; thus
$$
A=\sum_{j=2}^4 A_j \d x^j, \qquad \phi =\sum_{j=1}^3 \phi_j\,\d x^j,
$$ 
where the coefficient matrices $A_j$ and $\phi_j$ are functions of $x^2,x^3,x^4$ only.  With these assumptions, the KW equations take a nice form 
efficiently described using the three operators
\begin{equation}
\label{doff}
\begin{aligned}
\D_1 & = D_2+iD_3=\frac{\partial}{\partial x^2}+i\frac{\partial}{\partial x^3}  +[A_2+iA_3,\,\cdot\, ] \\
\D_2& =D_4-i[\phi_1,\,\cdot\,]= \frac{\partial}{\partial x^4} +[A_4-i\phi_1,\,\cdot\,],  \\
\D_3&=[\phi_2-i\phi_3,\,\cdot\,] , 
\end{aligned}  
\end{equation}
and the moment map
\begin{equation}
\label{zoff}
\mu=F_{23}-[\phi_2,\phi_3]-D_4\phi_1.
\end{equation}   
The KW equations with $A_1=\phi_4=0$ and all fields independent of $x^1$ then reduce to a set of ``complex equations'' 
\begin{equation}
\label{woff} 
[\D_i,\D_j]=0,~~1\leq i<j\leq 3,
\end{equation}
and a moment map condition
\begin{equation}
\label{old} 
\mu := \frac{i}{2} \sum_{j=1}^3 [\D_j, \D_j^\dag] = 0.
\end{equation}  
This is called the reduced system or the reduced KW equations. 

The complex equations are invariant under complex-valued gauge transformations 
\begin{equation}\label{zolf}
\D_i\to g\D_i g^{-1}, 
\end{equation}
where $g$ is a map from $\R^4_+$ to the complexified Lie group $G_\C$, while the moment map condition is invariant only under $G$-valued gauge transformations.
With suitable boundary conditions, one can understand the solutions of the complex equation, modulo $G_\C$-valued gauge transformations.  As in many somewhat
similar problems, one hopes to prove that 
solutions of the reduced system modulo $G$-valued gauge transformations are in bijective correspondence with solutions of
the complex equations modulo $G_\C$-valued gauge transformations, a much simpler problem. 

Now set $z=x^2+i x^3$,  $r=|z|$, $y=x^4$, $\uphi=\phi_2-i\phi_3$, so in particular $\D_3= [\uphi, \cdot]$. Using a complex gauge transformation, we 
may assume that 
\begin{equation}\label{murky}
\D^{(0)}_1 =\frac{\partial}{\partial\bar z}, \quad \D^{(0)}_2 =\frac{\partial}{\partial y}, \quad \D^{(0)}_3=\uphi_0(z), 
\end{equation}
where $\uphi_0(z)$ is holomorphic in $z$ and independent of $y$.  
We now take 
\begin{equation}\label{urky}
\uphi_0(z)=z^\r\begin{pmatrix} 0 & 1\cr 0 & 0 \end{pmatrix}
\end{equation}
Here $\r$ is a nonnegative integer that should be interpreted as a dominant weight of the Langlands
or GNO dual group, which we will call $G^\vee$.  
Thus, $\uphi_0(z)$ is a regular nilpotent element of the complex Lie algebra $\frak{sl}(2,\C)$ except at $z=0$,
where it vanishes.   There is a subtlety here as the dual of $SU(2)$ is  $SO(3)$, and vice-versa.  If $\r$ is even,
it can be interpreted as a weight of $G^\vee=SO(3)$, and this means that one gets a good model solution for
$G=SU(2)$.  But if $\r$ is odd, we have to interpret it   as a weight of $G^\vee=SU(2)$, meaning that the model
solution is well-behaved only for $G=SO(3)$.  In practice, it is useful to describe the model solution in the language
of $SU(2)$, and then explain why for $\r$ odd, the formulas should actually 
be interpreted in terms of $SO(3)$ gauge theory.  

We then solve for the $SL(2,\C)$-valued gauge transformation $g$ such that $\D_i=g \D_i^{(0)} g^{-1}$ satisfy the moment map condition. This
leads to the expression 
\begin{equation}\label{turfew}g=\begin{pmatrix}e^{v/2}& 0 \cr
0 & e^{-v/2}\end{pmatrix},~~~ e^v=\frac{2(\r+1)  }{(\sqrt{r^2+y^2}+y)^{\r+1}-(\sqrt{r^2+y^2}-y)^{\r+1}}.
\end{equation} 
The solution then takes the explicit form
\begin{equation}\label{expsol}
\begin{aligned}
A&=  -r\partial_rv(r,y)\,\d\theta\begin{pmatrix}\frac{i}{2}&0\cr 0&-\frac{i}{2}\end{pmatrix}\cr\uphi&=e^v z^\r\begin{pmatrix}0&1\cr 0&0\end{pmatrix}=  \frac{2(\r+1)z^\r  }{(\sqrt{r^2+y^2}+y)^{\r+1}-(\sqrt{r^2+y^2}-y)^{\r+1}}\begin{pmatrix}0&1\cr 0&0\end{pmatrix}
\cr\phi_1&=  -\partial_yv(r,y)\,\begin{pmatrix}\frac{i}{2}&0\cr 0&-\frac{i}{2}\end{pmatrix},
\end{aligned}
\end{equation}
with all other fields vanishing.  In spherical coordinates 
\begin{equation}\label{polar}
x^4=\rho\cos\psi, ~x^2=\rho\sin\psi\cos \theta,~x^3=\rho\sin\psi\sin\theta,
\end{equation}
this may be expressed as
\begin{equation}
\label{more}
\begin{aligned}A&=-(\r+1){\sin^2\psi}\frac{\left(1+\cos\psi\right)^{\r}-\left(1-\cos\psi\right)^{\r}}{\left(1+\cos\psi\right)^{\r+1}-\left(1-\cos\psi\right)^{\r+1}}\,\d\theta\begin{pmatrix}\frac{i}{2}&0\cr 0&-\frac{i}{2}\end{pmatrix}\cr
\uphi&=\frac{2(\r+1)}{\rho}\frac{\sin^r\negthinspace\psi \,e^{i\r\theta}}{\left(1+\cos\psi\right)^{\r+1}-\left(1-\cos\psi\right)^{\r+1}
}\begin{pmatrix}0&1\cr 0&0\end{pmatrix}\cr \phi_1&= -\frac{\r+1}{\rho}\frac{\left(1+\cos\psi\right)^{\r+1}+\left(1-\cos\psi\right)^{\r+1}}
{\left(1+\cos\psi\right)^{\r+1}-\left(1-\cos\psi\right)^{\r+1}}\begin{pmatrix}
\frac{i}{2}&0\cr 0&-\frac{i}{2}\end{pmatrix} \end{aligned}
\end{equation}
Scale-invariance means that $\uphi$ and $\phi_1$ are $1/\rho$ times functions of $\psi$ only (so that the corresponding one-forms $\uphi\d z$ and $\phi_1\d y$ 
are invariant under scaling), while $A_\theta$ is a function of $\psi$ alone. Rotation-invariance means that a constant rotation  $\theta\to \theta+c$ 
can be compensated by a constant diagonal gauge transformation $\mathrm{diag} (e^{-ic\r/2},e^{ic\r/2})$.

To see that this solution obeys the standard Nahm pole boundary condition away from $z=0$, and also
to see the role of the parity of $\r$, we observe that for  $z\not=0$, one has
\begin{equation}\label{nexp}\uphi\sim\left(\frac{z}{\bar z}\right)^{\r/2} \frac{1}{y}  \begin{pmatrix}0&1\cr 0&0\end{pmatrix},~~~y\to 0.
\end{equation}
This is converted to 
\begin{equation}\label{wexp} 
\uphi\sim \frac{1}{y}  \begin{pmatrix}0&1\cr 0&0\end{pmatrix},~~~y\to 0 ,
\end{equation}
by the gauge transformation $h=\mathrm{diag}((z/\bar z)^{-\r/4}, (z/\bar z)^{\r/4})$; furthermore, $h$ also makes $A=0$ and $\phi_1 \sim\mathrm{diag}(i/2y,-i/2y)$ 
as $y\to 0$, which gives the rest of the standard Nahm pole solution.   Observe that if $\r$ is odd, $h$ is only single-valued as an $SO(3)$,
but not an $SU(2)$, gauge transformation, so in that case the model solution \eqref{expsol} is gauge-equivalent to the standard Nahm pole solution 
along the boundary away from $z=0$ only for the gauge group $SO(3)$, while if $\r$ is even, then $h$ is single-valued and hence \eqref{expsol} is 
gauge-equivalent to this standard solution away from $z=0$ for either gauge group $SU(2)$ or $SO(3)$. 

The generalized Nahm pole boundary condition for a field $(A,\phi)$ on a manifold with boundary $X$ with knot $K \subset W = \del X$
can now be defined, at least informally; a quantitative definition will be given later after the calculation in Section 5 of the indicial roots of the problem.
The idea is simply that the model solution above can be transplanted to each fiber of $N^+K$. These fibers
are three-dimensional half-spaces with distinguished point on the boundary, the origin.  Transplanting the solutions to $N^+K$  is well-defined independent of the framing of $K$ 
by virtue of the rotational invariance of the model solutions about the vertical axis in each half-space.  This defines an approximate solution 
along the entire knot, which we call $(A^K, \phi^K)$.  A field $(A,\phi)$ on $X$ is then said to satisfy the generalized Nahm pole  boundary
condition at $K$ provided
\begin{equation}\label{gnpbc}
(A,\phi) \sim (A^K, \phi^K)\ \ \mbox{at}\ K.
\end{equation}
The clarification of \eqref{gnpbc} given later involves specifying the rate of convergence in this formula. This boundary condition along $K$ is supplemented by the 
Nahm pole boundary condition \eqref{npbc} at $W \setminus K$, as described fully in \cite{MW}.  

\section{The Linearized KW operator}\label{linearized}
We next describe the structure of the linearization $\LKW$ of the KW equations around the knot $K$. 

Let $(X,g)$ be a four dimensional Riemannian manifold with boundary $W$, $E$ an $\mathrm{SU(2)}$ or $\mathrm{SO(3)}$ bundle over $X$,
and $K \subset W$ a closed knot or link.  We now use local coordinates $\vec x = (x^1, x^2, x^3, x^4)$ near any point of $K$, where $x^1 = t$ is arclength 
along $K$, $x^2, x^3$ are Fermi coordinates around $K$ for the restriction of $g$ to $W$ around $K$ in $W$, and $y = x^4$ is geodesic distance 
from $W$, all with respect to $g$.  It is often more convenient to use the corresponding cylindrical coordinates $\rho = |(y, x^2, x^3)|$ and 
$\omega = (\omega_0, \omega') = (y, x^2, x^3)/\rho$ in the hemisphere $S^2_+$. We also use spherical coordinates \eqref{polar}, but
replacing $\psi$ by $s = \pi/2 - \psi$ to emphasize the fact that the important singularity is at $s=0$. Thus 
$\omega = (\sin s, \cos s \cos \theta, \cos s \sin \theta)$. In these coordinates
\begin{equation}
g = d\rho^2 + \rho^2( \d s^2 + \cos^2 s \d \theta^2) + \d t^2 + \ \mbox{higher order terms,} 
\label{g}
\end{equation}
where the remainder corresponds to the higher order terms in $\rho$ in the Taylor expansion of $g$. 

As usual, the linearization at a solution $(A,\phi)$ is obtained by placing $A+\varepsilon a ,\phi+\varepsilon \varphi$ into the KW equations, expanding 
in $\varepsilon$, and throwing away all terms of order $\varepsilon^2$ and higher. We shall supplement $\KW$, and hence $\LKW$, with the gauge 
condition 
 \begin{equation}\label{suppl}
S=\sum_{i=1}^4\left(D_ia_i+[\phi_i,\varphi_i]\right)=\sum_{i=1}^4\left([\partial_i+A_i,a_i] +[\phi_i,\varphi_i]\right) = 0,
\end{equation}
or, in invariant form, 
\begin{equation}\label{suppl2}
d_A^* a + \star [ \phi, \star \varphi] = 0
\end{equation}
to obtain an elliptic operator. Henceforth the ``linearized KW equations'' refer to the linearization of the KW equations supplemented 
by this gauge condition. 

Now suppose that $(A,\phi)$ is a solution of the gauged KW equations on $(M,K)$, with a Nahm pole singularity along $W$ 
away from $K$ and with the generalized Nahm pole condition along $K$.  
Recall from \cite{MW} that $\LKW$ takes a simpler form when we identify $\varphi$ with its Hodge dual $\star \varphi$. This identifies 
$\LKW$ with the twisted Hodge-de Rham operator
$$
\d + \d^*:  \calC^\infty(M, (\Lambda^1 \oplus \Lambda^3) \otimes \ad \, \mathfrak g) \longrightarrow
\calC^\infty(M, (\Lambda^0 \oplus \Lambda^2 \oplus \Lambda^4) \otimes \ad \, \mathfrak g),
$$
up to lower order terms, so in total
\begin{equation}
\begin{split}
\LKW(a, \star \varphi) = & \big(\d_A^* a + \star [ \phi,  \star \varphi], \ \  \d_A a + [\phi, \star (\star \varphi)] + 
\d_A^* \star \varphi + \star [a, \phi], \ \  \d_A \star \varphi + [a, \star \phi]\big) \\[0.5ex]
&  = \big(\d_A + \d_A^*\big) \big(a, \star \varphi\big) + 
\big( \star [ \phi, \star \varphi],  [\phi, \star (\star \varphi)] + \star [a, \phi], [a, \star \phi]\big).
\end{split}
\label{linLKW2}
\end{equation}
Near $W$ but away from $K$, 
$$
\LKW = \d + \d^* + \frac{1}{y} B_0,
$$
while in cylindrical coordinates near $K$, 
\begin{equation}
\LKW = \d + \d^* + \frac{1}{\rho s} B_0\, ;
\label{1LKW}
\end{equation}
here $B_0$ is an endomorphism determined by $g$, $A$ and $\phi$.   We remark finally that the lower order term in the
gauge condition \eqref{suppl2} is chosen so that $\LKW$ has a particularly nice form. 

\section{Indicial Operator}\label{indicial}
In this section we consider the first basic model for the linearized KW operator $\LKW$, called the indicial operator.  This is an infinitesimal
model for $\LKW$ at points of $W$.  The indicial operator at points of $W \setminus K$ was studied in \cite{MW}, and we begin with a brief review 
of that case.  For such points, the indicial operator is an ODE in the normal variable which is homogeneous of degree $-1$, and its solutions 
correspond to possible leading terms in formal expansions of solutions to the full equation $\LKW \Psi  = 0$.  The indicial operator at points of $K$
is more complicated in that it is a PDE in the directions normal to $K$, but its solutions play an identical role as possible leading terms
in expansions of solutions to $\LKW \Psi = 0$ near $K$.  There is another model for $\LKW$ which captures more delicate features and which
will be introduced in Sec.\ 7.  In any case, we now describe these indicial operators and then turn to the main task of computing the indicial roots of
the indicial operator at $K$. 

If $(\vec{x}, y)$ is a local coordinate system near $q \in W \setminus K$ with $y = x^4 \geq 0$, we say that $\lambda$ is an indicial root
for $\LKW$ at $q$ if there exists some $(a, \varphi)$ defined in a neighborhood of $q$ such that 
$$
\LKW ( y^\lambda (a, \varphi)) =  \mathcal O( y^{\lambda}).
$$
Note that since $\LKW$ is of order $1$, for general $\lambda$ one expects a right hand side which is $\calO(y^{\lambda-1})$, so we are really solving 
some sort of eigenvalue equation to make the coefficient of $y^{\lambda-1}$ vanish. This is an entirely algebraic problem, and the value of
$(a,\varphi)$ at $q$ solves a generalized eigenvector equation.  The nontrivial task of determining these eigenvectors and eigenvalues 
is treated in \cite{MW}. 

This indicial data corresponds to exact solutions of the indicial operator at $q$,
\begin{equation}\label{indop1}
I_q(\LKW) = U_y \del_y + \frac{1}{y} U_0
\end{equation}
where $U_y$ and $U_0$ are constant elements of $\mbox{End}(\Lambda^1 T^*_q M \otimes \ad (E))$. This is 
obtained by discarding all of the higher order terms in $\LKW$ as well as all terms with tangential ($\vec{x}$) derivatives; equivalently, 
it is precisely the linearized KW operator $\LKW$ in the model case $\RR^4_+$ with the Euclidean metric, and acting on fields 
independent of $\vec{x}$.  Comparing 
with the previous paragraph, we see that if $\Psi$ is smooth in $\vec{x}, y$ with $\Psi(\vec{x},0) = \Psi_0(\vec{x})$, then
$$
\LKW( y^\lambda\Psi) = I_q(\LKW) (y^\lambda \Psi_0) + \calO(y^\lambda),
$$
so $\lambda$ is an indicial root with corresponding eigenvector if and only if
$$
(U_y \lambda + U_0) \Psi_0 = 0. 
$$
For $G = {SU}(2)$ or ${SO}(3)$, the indicial roots are $-2$, $-1$, $1$, and $2$ \cite{MW}.  Notice that these are 
independent of the location of $q \in W$.   

The relevance of these values lies in a theorem from \cite{Ma}, stated in \cite{MW}, that if $\LKW (a,\varphi) = 0$ in some neighborhood of the origin,
and if $|(a, \varphi)| \leq C y^{-1 + \varepsilon}$ for some $\varepsilon > 0$, then $(a, \varphi)$ admits a polyhomogeneous expansion, i.e.,
an asymptotic expansion with tangentially smooth coefficients in increasing but possibly nonintegral powers of $y$, starting
from the first positive indicial root: 
$$
(a, \varphi) \sim  \sum_{\ell \geq 0} (a_\ell, \varphi_\ell) y^{1 + \ell}
$$
with coefficients $(a_\ell, \varphi_\ell)$ independent of $y$ and depending smoothly on the tangential variables $\vec{x}$.  (Since the positive indicial 
roots differ by an integer in our setting, there is also the possibility of terms of the form $y^{1 + j} (\log y)^\ell$ for $j, \ell \geq 1$, but these do not 
contribute in any significant way below, so we mostly omit them for simplicity of notation.) 

We refer to the condition that $(a, \varphi)$ is allowed to blow up no faster than $y^{-1 + \varepsilon}$ as the Nahm pole boundary condition
(when $K = \emptyset$). Analogous to a standard nondegenerate elliptic boundary condition, this guarantees that $(a,\varphi)$ decays like $y$ 
and is smooth (or at least polyhmogeneous) up to the boundary. 

Now let us turn to the indicial roots near points of $K$. Using the cylindrical coordinates $(\rho, t, \omega)$ introduced earlier, 
we say that $\lambda$ is an indicial root of $\LKW$ at $q \in K$ if there exists some $(a_0, \varphi_0)$ depending only on $\omega$ so that
$$
\LKW (\rho^\lambda (a, \varphi)) = \mathcal O (\rho^\lambda), 
$$
where $(a,\varphi)$ is any smooth (in $\rho$) extension of $(a_0, \varphi_0)$. As before, because $\LKW$ is to leading order homogeneous
of degree $-1$ with respect to dilation in $\rho$, we expect only that the right hand side is 
$\mathcal O( \rho^{\lambda -1 })$ for general $\lambda$, so the indicial equation is simply the coefficient 
of $\rho^{\lambda-1}$ here and the indicial data $\lambda$, $(a_0, \varphi_0)$ correspond to solutions which make this 
coefficient vanish.  This is again a generalized eigenvalue problem, but this time for an elliptic operator rather than a matrix.
Note that the terms involving differentiations with respect to $t$ do not contribute to this coefficient, 
nor do the terms arising from the Taylor expansions in $\rho$ of all the other coefficients except the leading term. Discarding
all of these terms leads to the indicial operator at $q \in K$, 
\begin{equation}\label{indop2}
I_q(\LKW) = B_\rho \del_\rho + \frac{1}{\rho} \calJ_S,
\end{equation}
where $B_\rho$ is a constant Clifford-multiplication endomorphism and $\calJ_S$ is a first order elliptic operator on the hemisphere $S^2_+$.
As in the case when $q \notin K$, this equals the linearized KW operator in the model case, and acting on $t$-independent fields. 
The rotation invariance of the model knot singularity implies that $\calJ_S$ is invariant with respect to rotations in $\theta$.  Furthermore,
$$
\LKW ( \rho^\lambda (a, \varphi)) = ( B_\rho + \lambda \calJ_S)(a_0, \varphi_0) \rho^{\lambda-1} + \calO(\rho^\lambda),
$$
whence the indicial equation at $ q\in K$, 
\begin{equation}\label{indeqn}
(B_\rho + \lambda \calJ_S) (a_0, \varphi_0) = 0. 
\end{equation}
Note that while the fields $(a_0, \varphi_0)$ depend only on $\theta$ and $\psi$, they may have $\d\rho$ components. 

Our goal is to describe the indicial operator at $q \in K$ and analyze its spectrum in sufficient detail to show that there are no indicial 
roots in the semi-open interval $[-1,0)$. In fact, we will actually see that only a rather restricted class of perturbations lead to indicial 
roots in the larger interval $[-2,1]$. This bound on the indicial roots will be important in the ensuing analysis.  

The computation of the indicial roots of $\LKW$ is carried out using two separate techniques.  The first 
uses the linearization of the Weitzenb\"ock  formula that will be given in section 5 below, which equates $\LKW^\dag \LKW$ 
with a certain Laplace-type operator $\hD$.  In the particular case when we linearize around the model
knot solution, this operator acts separately on $a_1$ and $\varphi_4$ and decouples from the other 
components, and we call this specialization $\hD_0$.  This decoupling uses that $A_1 = \phi_4 = 0$. If 
in addition the fields $(A,\phi)$ are independent of $x_1$, as holds for the model knot solution, then
$\hD_0$ reduces to an operator $\hD_{0,3}$ in $\RR^3_+$.  We shall calculate the indicial roots of $\hD_0$, or
equivalently, $\hD_{0,3}$; the indicial roots for $\LKW$ in these directions are amongst these, but 
we do not determine which are roots for $\hD_0$ only, but not $\LKW$. On the other hand, perturbations 
where the components $a_1, \varphi_4$ remain zero correspond to perturbations within the three-dimensional 
formulation of the KW equations \eqref{doff} involving the operators $\D_i$, $i=1, 2, 3$. For these our 
computation give the indicial roots of $\LKW$ in those directions directly. 

To summarize,  we consider two classes of perturbations: 

(1) Perturbations of Type I have either $a_1\not=0$ or $\varphi_4\not=0$. These are treated using 
the operator $\hD$. 

(2) Perturbations of Type II are those that preserve the condition $A_1=\phi_4=0$ and can be analyzed using \eqref{woff} 
and \eqref{old}. 

\subsection{Perturbations Of Type I}\label{typetwo}
Consider a solution $(a,\varphi)$ to the indicial equation $I_\lambda(\LKW)(a,\varphi) = 0$, where either $a_1 \neq 0$ or 
$\varphi_4 \neq 0$. 
Our analysis is based on the observation that $I_\lambda(\LKW)^\dagger I_\lambda(\LKW) (a,\varphi) = 0$, and this leads to a 
particularly simple equation for $a_1$ and $\varphi_4$.

We calculate $\hD$ using the action functional 
\begin{equation}\label{actfun}
I(A,\phi) =-\frac{1}{2}\int\d^4x \,\Tr \left(\V_{ij}\V^{ij}+(\V^0)^2\right)
\end{equation}
of the underlying supersymmetric gauge theory, where $\V_{ij}=F_{ij}- [\phi_i,\phi_j]+\epsilon_{ij}{}^{kl}D_k\phi_l$, 
$\V^0=D_i\phi^i$, and $\V_{ij}=\V^0=0$ are the KW equations. This action vanishes if and only if ${\KW}(A,\phi) = 0$.  
Now observe the following elementary fact: let $N$ be some nonlinear functional on fields $\Psi$ (for simplicity we 
assume that its nonlinear terms are only quadratic) and $I(\Psi) = \int |N(\Psi)|^2$ the corresponding action.  
If $N(\Psi_0) = 0$, then for any $\psi$, 
$$
N(\Psi_0 + \varepsilon \psi) = N(\Psi_0) + \varepsilon \LKW' \psi + \varepsilon^2 Q(\Psi_0, \psi)
$$ 
where $Q$ is quadratic in $\psi$ and $\LKW' = DN|_{\Psi_0}$, hence $I( \Psi_0 + \varepsilon \psi) =  
\varepsilon^2 \int \left|\LKW' \psi\right|^2 + \calO(\varepsilon^3)$. This discussion is not quite suitable 
in our setting because $I$ is invariant under gauge transformations, so $\LKW'$ is not elliptic. 
This can be remedied either by augmenting $\LKW'$ with the equation $S=0$, where $S$ is the 
gauge-fixing equation of $\Psi$ relative to $\Psi_0$, or else by adding $-\int d^4x \, \Tr S^2$ to $I(A,\phi)$ 
and then carrying out the calculations above. Either way, if $I_2(\psi)$ denotes the coefficient of $\varepsilon^2$ in
the expansion of $I$, then critical points of $I_2$ also satisfying the gauge condition are solutions of 
$\LKW^\dagger \LKW \psi = 0$, where $\LKW = \LKW' + S$ is the linearized gauged KW operator. 

Now, proceeding as in Eqns.\ (2.51-52) of \cite{MW}, and preferably using the formulas for $I(A,\phi)$ given 
below in eqns. \eqref{baction} or \eqref{zoffbo2} (as this simplifies the computations), one finds that
\begin{equation}\label{zite}
I_2=I_{2,0}+I_{2,1}+I_{2,2},
\end{equation}
where 
\begin{equation}
\label{ztwo}
\begin{aligned}I_{2,0}&=-\int \d^4x\,\Tr\sum_{i,j=1}^4\left((D_i a_j)^2+(D_i\varphi_j )^2+[\phi_i,a_j]^2+[\phi_i,\varphi_j]^2\right) \\
 I_{2,1} &=-2\int \d^4 x \,\Tr \sum_{i,j=1}^4\left(F_{ij}[a_i,a_j]+[\phi_i,\phi_j][\varphi_i,\varphi_j]+2D_i\phi_j [a_i,\varphi_j] \right) \\ 
I_{2,2}&= \int \d^4x\, \Tr\,S^2.
\end{aligned}
\end{equation}
Boundary terms have been dropped since they do not contribute to the Euler-Lagrange equations.  Restricting to gauge-fixed variations, 
we may also drop $I_{2,2}$. 
We now define $\hD$ to be the second order elliptic operator obtained as the Euler-Lagrange equation for $I_{2,0} + I_{2,1}$ at
the model knot solution. 

If we expand around a solution of the KW equations for which $\phi_4 \equiv 0$, then $\varphi_4$ is absent in $I_{2,1}$, so in fact
the Euler-Lagrange equation for $\varphi_4$ can be derived from $I_{2,0}$ alone, and is calculated to be
\begin{equation}\label{iffy}
-\left(\sum_{i=1}^4 D_i^2+\sum_{j=1}^3[\phi_j,[\phi_j,\cdot]]\right)\varphi_4=0,~~~
D_i=\frac{\partial}{\partial x^i}+[A_i,\cdot]. 
\end{equation} 

In a similar way, if we expand around a solution of the KW equations for which $F_{1i}$ and $D_1\phi_i$ vanish identically for $i=2,3,4$ (this is the case, in particular, if the solution 
is invariant under translations in the $x^1$ direction and also has $A_1=0$), then $a_1$ does not appear in $I_{2,1}$, so the Euler-Lagrange equation 
for  $a_1$  is again derived just from $I_{2,0}$ and is 
\begin{equation}\label{jiffy} 
-\left(\sum_{i=1}^4 D_i^2+\sum_{j=1}^4[\phi_j,[\phi_j,\cdot]]\right)a_1=0. 
\end{equation}
Note that these two equations coincide if all conditions are satisfied ($\phi_4=F_{1i}=D_1\phi_i=0$).  
 
The model solution for a knot has the properties assumed in the last two paragraphs, and we shall denote as   $\hD_0$ the common operator
appearing in \eqref{iffy} and \eqref{jiffy} in expanding around that solution. It is then straightforward to describe the indicial 
equation for Type I perturbations around this solution. 
As explained earlier, we restrict to perturbations $a,\varphi$ which are independent of $x^1$ and obey the three-dimensional reduction
\begin{equation}\label{niffy} 
\hD_{0,3} \Psi=0, \quad \mbox{where}\qquad  \hD_{0,3} =-\left(\sum_{i=2}^4 D_i^2+\sum_{j=1}^3[\phi_j,[\phi_j,\cdot]]\right). 
\end{equation}
Here $\Psi$ is either $a_1$ or $\varphi_4$.  The operator $\hD_{0,3}$ appears frequently below. 

In polar coordinates \eqref{polar}, the ordinary Laplacian $\Delta=-\sum_{i=2}^4\partial_i^2$ takes the form 
\begin{equation}\label{liffy}
\Delta=-\left(\frac{\partial^2}{\partial\rho^2}+\frac{2}{\rho}\frac{\partial}{\partial\rho}
+\frac{1}{\rho^2}\left(\frac{\partial^2}{\partial\psi^2}+\frac{\cos\psi}{\sin\psi}\frac{\partial}{\partial\psi}+
\frac{1}{\sin^2\psi}\frac{\partial^2}{\partial\theta^2}\right)\right). 
\end{equation}
In the presence of the gauge field $A$ , we must simply replace derivatives by covariant derivatives.   For the model solution for a knot, 
this is particularly simple because by \eqref{expsol}, $A$ is a multiple of $\d\theta$, so defining
\begin{equation}\label{wiffy}
D_\theta=\partial_\theta+[A_\theta,\cdot],~~~ A_\theta=-r\partial_rv \begin{pmatrix}
\frac{i}{2}& 0 \cr 0 & -\frac{i}{2}\end{pmatrix}, 
\end{equation}
where $v$ is the function in \eqref{expsol}, we have
\begin{equation}\label{uplift}
-\sum_{j=2}^4 D_j^2=-\left(\frac{\partial^2}{\partial\rho^2}+\frac{2}{\rho}\frac{\partial}{\partial\rho}
+\frac{1}{\rho^2}\left(\frac{\partial^2}{\partial\psi^2}+\frac{\cos\psi}{\sin\psi}\frac{\partial}{\partial \psi}+\frac{D_\theta^2}{\sin^2\psi}\right)\right). 
\end{equation}

The other term in \eqref{niffy} is 
\begin{equation}\label{omolt}
\Delta_\phi=-\sum_{i=1}^3[\phi_i,[\phi_i,\cdot]] =-[\phi_1,[\phi_1,\cdot]]+\frac{1}{2}\left([\uphi^\dagger,[\uphi,\cdot]]+[\uphi,
 [\uphi^\dagger,\cdot]]\right)=\frac{N_S}{\rho^2},
\end{equation}  
where $N_S$ is an angle-dependent linear transformation of the real Lie algebra $\frak{su}(2)$ and $\uphi$ is as in \eqref{more}. 
Note that $N_S$ is strictly positive.
One sees directly from \eqref{more} that  $N_S$ depends only on $\psi$. Thus 
\begin{equation}\label{telox}
\hD_{0,3} =-\frac{\partial^2}{\partial\rho^2}-\frac{2}{\rho}\frac{\partial}{\partial\rho} +\frac{M_S}{\rho^2},
\end{equation}
where 
\begin{equation}\label{welox}
M_S=-\frac{\partial^2\, }{\partial\psi^2}-\frac{\cos\psi}{\sin\psi}\frac{\partial}{\partial\psi}-\frac{D_\theta^2}{\sin^2\psi}+N_S
\end{equation}
is an operator acting on the hemisphere $S^2_+$ which is invariant with respect to rotations in $\theta$. The apparent singularity at $\psi = 0$ is only a 
polar coordinate singularity, but on the other hand, because of the generalized Nahm singularity in $\phi_1$ and $\uphi$, the transformation $N_S$ blows 
up like $1/\cos^2 \psi$ at the boundary of this hemisphere.  

We shall need to solve equations of the form $M_S \Psi_0 = f$ and $M_S \Psi_0 = \lambda \Psi_0$ below, so we review some analytic properties
of this operator.  To emphasize the singularity at $\psi = \pi/2$, let us revert to the coordinate $s = \pi/2 - \psi$ and write
$$
M_S = -\frac{\partial^2\, }{\partial s^2} + \frac{\sin s}{\cos s}\frac{\partial\, }{\partial s}-\frac{D_\theta^2}{\cos^2 s}+N_S
$$
Using \eqref{omolt} and the asymptotic expressions for the $\phi_j$ (or equivalently, $\phi_1$ and $\uphi$) we compute that 
$$
N_S \sim \frac{2 \, \mbox{Id}}{s^2} \ \ \mbox{as}\ \ s \searrow 0,
$$
and hence 
$$
M_S = -\frac{\del^2\,}{\del s^2} - D_\theta^2 + \frac{2 \, \mbox{Id}}{s^2} + \calO(s). 
$$

Recall from \cite{Ma} and \cite{MW} that $s^2 M_S$ is an elliptic uniformly degenerate operator. Much of its behavior is dictated by 
its own indicial roots, i.e., the values $\sigma$ such that $M_S ( s^\sigma) = \calO( s^{\sigma-1})$
(the expected order of decay is $s^{\sigma-2}$ since $M_S$ is second order, so this represents the same sort of leading
order cancellation as we have been discussing). A brief calculation shows that the indicial roots are $-1$ and $2$. We now collect 
a set of results from \cite{Ma} regarding the mapping properties of such operators on weighted $L^2$ and H\"older spaces: 
\begin{prop}\label{MSprop}
The operators
$$
M_S:  s^\nu H^2_0( S^2_+; \d s \, \d\theta) \longrightarrow s^{\nu-2} L^2( S^2_+; \d s \, \d\theta)
$$
$$
M_S:  s^\mu \calC^{2,\alpha}_0( S^2_+) \longrightarrow s^{\mu-2} \calC^{0,\alpha}_0( S^2_+)
$$
are invertible provided $-1/2 < \nu < 5/2$, $-1 < \mu < 2$. There is a unique self-adjoint realization of 
$$
M_S:  L^2 (S^2_+; \d s \, \d\theta) \longrightarrow L^2(S^2_+; \d s \, \d\theta),
$$
which has domain $s^2 H^2_0(S^2_+; \d s \, \d\theta)$; it has discrete spectrum.  If $f \in s^{\mu-2} \calC^\infty(S^2_+)$ 
for $-1 < \mu < 2$ and $(M_S + V) \Psi_0 = f$, where $\Psi_0 \in s^{-1/2 + \varepsilon} L^2$ for some $\varepsilon > 0$, 
then $\Psi_0 \in s^\mu \calC^\infty(S^2_+)$.  In particular, $\Psi_0(0,\theta) = 0$ if $\mu > 0$. 
\end{prop}

We comment only briefly on the proof. The techniques of \cite{Ma} allow one to construct a `parametrix', an approximate
inverse to any of these mappings (we discuss the parametrix technique at some length in Sec.\ 7 below). This leads directly
to the Fredholm properties, the regularity estimates and the essential self-adjointness. In particular, solutions of
$M_S \Psi_0 = \lambda \Psi_0$ or $M_s \Psi_0 = f$ with $f \in \calC^\infty$ admit polyhomogeneous expansions as $s \to 0$
with leading term $s^\sigma$ where $\sigma$ is the smallest indicial root in $L^2$, which in this case is $1$.  The salient 
conclusion for us is that solutions must vanish at $s=0$. The fact that $M_S$ is not just Fredholm but actually invertible
follows from its symmetry and integration by parts to eliminate its nullspace.

We can find fields $\Psi$ solving $\hD_{0,3} \Psi = 0$ and which are homogeneous in $\rho$, i.e., $\Psi(\rho,\psi,\theta)=\rho^\lambda 
\Psi_0(\psi,\theta)$, by calculating 
\begin{equation}\label{zelox}
\hD_{0,3} (\rho^\lambda \Psi_0) = 0  \Longleftrightarrow (M_S - \lambda(\lambda+1)) \Psi_0 = 0.
\end{equation}
Thus we must choose $\Psi_0$ to be an eigenvector for $M_S$ with eigenvalue $\gamma = \lambda(\lambda + 1)$, or equivalently, 
\begin{equation}\label{weloxe}
\lambda=-\frac{1}{2}\pm \sqrt{\gamma+\frac{1}{4}}. 
\end{equation}

Any Type I indicial root has the form \eqref{welox}, for some choice of sign of the square root, and some eigenvalue $\gamma$ of $M$. 
To prove the converse,  that every number $\lambda$ arising this way is an indicial root of the underlying KW problem, one must show that 
every $a_1$ or $\varphi_4$ that is annihilated by $\hD_{0,3}$ can be extended to a full set of fields $a_i,\varphi_j$ obeying  the indicial 
operator of the linearized KW equations. We expect this to be true but have not shown it. 

It will be useful to know that there are no indicial roots in as large an interval around $-1/2$ as possible. For this, we must demonstrate
a lower bound for the smallest eigenvalue of $M_S$.  Note that this smallest eigenvalue decreases if we drop the nonnegative terms 
$-D_\theta^2/\sin^2\psi$ and $N_S(\psi)$, so in other words we shall compute the smallest eigenvalue of 
$M_0=-\partial_\psi^2-\cot\psi \partial_\psi$.  This is nothing more than the scalar Laplacian on the hemisphere acting on 
$\theta$-independent functions, and its smallest $L^2$ eigenvalue is $2$, with corresponding eigenfunction $\cos\psi$.  
Hence the smallest eigenvalue of $M_S$ is strictly larger than $2$, which means, according to \eqref{weloxe}, 
that there are no indicial roots of Type I in the closed interval $[-2,1]$.  
 
It is likely possible to improve this estimate, possibly even by finding the eigenvalues of $M_S$ in closed form,  but  this lower bound
cannot be improved very much in the sense that the smallest eigenvalue actually does converge to $2$ as the parameter $\r\to\infty$.   Indeed,
consider the case that the wavefunction $\Psi_0$ is diagonal in the basis used in \eqref{more} and depends only on $\psi$. Then the terms in
$M_S$ involving $D_\theta$ and $\phi_1$ do not contribute, so the difference between $M_S$ and $M_0$ comes entirely from the terms
$[\uphi,[\uphi^\dagger,\cdot]]+[\uphi^\dagger,[\uphi,\cdot]]$ in $N_S$. However, $\uphi$ vanishes rapidly as $\r\to\infty$ except very near 
$\psi=\pi/2$, where the wavefunction $\Psi_0$  vanishes. Exploiting this, one can argue that for large $\r$, the smallest eigenvalue of $M_S$ 
for diagonal perturbations is $2+\mathcal{O}(1/\r)$. 
On the other hand, if $\Psi_0$ is upper- or lower-triangular, the smallest eigenvalue of $M_S$ grows with $\r$ and hence our bound on the  Type I
indicial roots is not very tight.

\medskip

The mapping properties of the operators $\hD_{0,3}$ may now be analyzed using the eigenfunction decomposition of $M_S$;
a similar but more complicated analysis also leads to the mapping properties of $\hD_0 = -D_1^2 + \hD_{0,3}$. 
\begin{prop}\label{hD0prop}
The operators
\begin{eqnarray*}
\hD_{0,3}: & \rho^\delta s^\nu H^2_{\ie}( \RR^3_+;  \rho^2 \d\rho \, \d \psi \, \d\theta) \longrightarrow \rho^{\delta-2} 
s^{\nu-2} L^2( \RR^3_+;  \rho^2 \d \rho \, \d \psi \, \d\theta)\\ 
\hD_0: & \rho^\delta s^\nu H^2_\ie( \RR^4_+;  \rho^2 \d\rho \, \d t \, \d \psi \, \d\theta) \longrightarrow \rho^{\delta-2} 
s^{\nu-2} L^2( \RR^4_+;  \rho^2 \d \rho \, \d t \, \d \psi \,\d\theta)
\end{eqnarray*}
and
\begin{eqnarray*}
\hD_{0,3}:  & \rho^\eta s^\mu \calC^{2,\alpha}_\ie( [\RR^3_+; \{0\}]) \longrightarrow \rho^{\eta -2} s^{\mu-2} \calC^{0,\alpha}_\ie( [\RR^3_+; \{0\}]) \\
\hD_{0}:  & \rho^{\eta} s^\mu \calC^{2,\alpha}_\ie( [\RR^4_+; \RR ]) \longrightarrow \rho^{\eta-2} s^{\mu-2} \calC^{0,\alpha}_\ie( [\RR^4_+; \RR])
\end{eqnarray*}
are invertible provided 
\begin{equation*}
\begin{split}
1/2 -\sqrt{1+\gamma_0} := \delta^- & < \delta < \delta^+ := 1/2 + \sqrt{1+\gamma_0},\  \\ 
-1-\sqrt{1+\gamma_0} := \eta^- & < \eta <  \eta^+ := -1 + \sqrt{1+\gamma_0}, 
\end{split}
\end{equation*}
where $\gamma_0$ is the smallest eigenvalue of $M_S$, and 
$$
-1/2 < \nu < 5/2,\ \   -1 < \mu < 2. 
$$
\end{prop}
The proof requires the analytic techniques developed in sections 8 and 9. We use notation and ideas which will be 
explained more carefully in those sections. The spaces $[\RR^3_+; \{0\}]$ and $[\RR^4_+; \RR]$ appearing in the H\"older 
spaces are the blowups of $\RR^3_+$ around the origin and $\RR^4_+$ around $\RR \subset \del \RR^4_+$, respectively; 
see section 8.1.  The function spaces decorated with a subscript $\ie$ are iterated edge spaces, and are defined in 
sections 9.2.2 and 9.3.6.  The ranges of weight parameters here are optimal and are based on the indicial root computations
in the present section. Sobolev and H\"older mapping properties are proved using the boundedness of parametrices for 
$\hD_{0,3}$ and $\hD_0$ as constructed in section 9.  These parametrices also lead to sharp regularity statements. Thus, if 
$f \in \rho^{\eta -2} s^{\mu-2} \calC^\infty([\RR^3_+; \{0\})$ for $\eta$ and $\mu$ in the ranges above,
and $\hD_{0,3}\Psi_0 = f$ where $\Psi_0 \in \rho^{\delta^- + \varepsilon} s^{-1/2 + \varepsilon} L^2$ 
for some $\varepsilon > 0$, then $\Psi_0$ is polyhomogeneous as $\rho \to 0$ with index set corresponding 
to the indicial root set \eqref{weloxe}, and smooth up to the boundary $s=0$ (vanishing there if $\mu > 0$).
Expansions of this type are explained more carefully in section 9. We have stated this result here because 
in the remaining calculations of indicial roots it is necessary to solve the equation $\hD_{0,3} \Psi_0 = f$ 
for specific right hand sides $f$, and this result shows that this is possible and that the solutions have
the same regularity as $f$. 
 
 \subsection{Perturbations Of Type II}\label{typeone}
We next consider the perturbations of Type II, i.e., perturbations which preserve the condition $A_1=\phi_4=0$. 
These can be described as perturbations of the three operators $\D_i$ that preserve the complex equation  $[\D_i,\D_j]=0$,
the moment map condition $\mu=\sum_i[\D_i,\D_i^\dagger] = 0$, and the gauge condition $S=0$.
 
\subsubsection{The Generic Case}\label{generic}
The obvious perturbations of the $\D_i$ which preserves their commutativity are given by conjugation, i.e., 
 \begin{equation}\label{forf}
\exp(-\varepsilon \chi) \D_i \exp(\varepsilon \chi)  =  \D_i + \varepsilon [\D_i,\chi] + \calO(\varepsilon^2),
\end{equation}
where $\chi$ is the generator of a complex gauge transformation, i.e., is valued in the complexification $\frak g_\C$ of the Lie algebra.
Most Type II deformations are of this kind and we call them roots of generic Type II.  We account for a few discrepancies later. 
  
The deformations $\chi$ have a somewhat different flavor depending on whether $\chi$ is valued in the real Lie algebra $\frak g$ or 
in $i\frak g$.  (When $\frak g$ is the space of traceless antihermitian matrices, these two cases correspond to $\chi$ being 
skew-hermitian or hermitian.)  If $\chi$ is real, it generates a $G$-valued gauge transformation that trivially preserves the 
moment map condition $\mu=0$ as well as the commutativity of the $\D_i$, but it might not preserve the gauge condition. 
Indeed, the equation $S=0$ becomes a second order differential equation for $\chi$.  On the other hand, if $\chi$ is imaginary, 
then it generates a $G_\C$-valued gauge transformation which preserves commutativity of the $\D_i$ and, perhaps unexpectedly, 
also the condition $S=0$. It is not a symmetry of the moment map, and the condition $\mu=0$ becomes a second order differential  
equation for $\chi$.  Curiously, the second order differential equation that comes from $S=0$ when $\chi$ is real coincides 
with the equation that comes from $\mu=0$ when $\chi$ is imaginary, and is the equation \eqref{niffy} which already appeared 
in the study of Type I deformations! 
 
First assume that $\chi$ is real.  Then the deformation \eqref{forf} corresponds to the fluctuations $a_i=D_i\chi$, $i=2,3,4$ and $\varphi_i=
[\phi_i,\chi]$, $i=1,2,3$.   The condition $S=0$ is $D_i a_i+[\phi_i,\varphi_i]=0$, which translates to
\begin{equation}\label{zeffo}
\hD_{0,3} \chi = -\left(\sum_{i=2}^4 D_i^2+ \sum_{i=1}^3[\phi_i,[\phi_i,\cdot]]\right)\chi=0,
\end{equation}
precisely the same operator as before!  

The mapping properties of $\hD_{0,3}$ in Proposition~\ref{hD0prop} show that $S=0$ is a good gauge condition.
The Nahm pole at a general boundary point is irreducible and thus not invariant under a gauge transformation which
is nontrivial at the boundary. Thus we consider only generators $\chi$ which vanish there.
We shall define the Nahm pole boundary condition so that in perturbing around any solution $\Psi=(A,\phi)$ which 
has the Nahm pole singularity at a general boundary point and a more subtle singularity such as we have described along
$K$, then the allowed perturbations $(a,\varphi)$ are less singular than $(A,\phi)$ both along the boundary and along $K$. 
This condition is only compatible with gauge transformations with generators $\chi$ vanishing along the boundary.  
Now, Proposition~\ref{hD0prop} implies that, acting on a space of generators $\chi$ which vanish to some small positive
order at the boundary, $\hD_{0,3}$ is invertible. Therefore, if $(a,\varphi)$ is a linear perturbation that satisfies the Nahm 
pole boundary conditions, then we can transform uniquely by a linearized gauge transformation 
$(a,\varphi)\to (a,\varphi)+ (-\d_A\chi, [\chi,\phi]),$ where $\chi$ vanishes on the boundary, so that $S=0$.  
Hence in the function space in which we  work here, it is equivalent either to define the linearized KW equations 
without a gauge condition and divide by gauge transformations that are trivial on the boundary, or to include $S=0$ as 
part of the linearized KW equations. We take the latter route in this paper.

Now consider the case that $\chi$ is imaginary.   Then the deformation \eqref{forf} amounts to
\begin{equation}
\label{sell}
\begin{aligned}
a_2 & = iD_3\chi, \quad a_3 = -iD_2\chi, ~~~ \quad a_4 =-i[\phi_1,\chi] \\ 
\varphi_1 & =iD_4\chi, \quad \varphi_2=-i[\phi_3,\chi], \quad \varphi_3 =i[\phi_2,\chi].
\end{aligned}
\end{equation}
For this variation, the gauge condition $S=0$ is trivially satisfied since the equations \eqref{sell} imply that $\sum_iD_i a_i+
\sum_i [\phi_i, \varphi_i]=i[\mu,\chi]$, which vanishes when $\mu=0$. Using \eqref{sell} again, we see that the variation in $\mu$ is 
\begin{equation}\label{merox}
\delta\mu=-\left(\sum_i D_i^2+\sum_i[\phi_i,[\phi_i,\cdot]]\right)\chi,
\end{equation}
so imposing that $\delta\mu=0$ leads once again to $\hD_{0,3} \chi=0$.    

We summarize and slightly expand on these facts as follows.  Consider any perturbation $\D_i \to \D_i + \delta \D_i$ corresponding to 
$D_i \to D_i + a_i$, $\phi_i \to \phi_i + \varphi_i$. Then by a brief calculation, 
\begin{equation}\label{linSmu0}
\sum_{j=1}^3 [ \D_j, \delta \D_j^\dag] = -S(a,\varphi) + i \delta \mu (a, \varphi).
\end{equation}
The paragraphs above show that if the perturbation comes from a real (skew-Hermitian) gauge transformation $\chi$, then 
\begin{equation}\label{linSmu1}
(-S + i\delta \mu)(a,\varphi) = \hD_{0,3} (\chi)
\end{equation}
while if $\chi$ is Hermitian, then
\begin{equation}\label{linSmu2}
(-S + i\delta \mu)(a,\varphi) = \hD_{0,3}(\chi) + i [\mu(a,\varphi), \chi]. 
\end{equation}
The perturbations by real or imaginary $\chi$  are related by the complex rotation 
$$
(a_2, a_3, a_4, \varphi_1, \varphi_2, \varphi_3) \mapsto ( i a_3, -i a_2, -i\varphi_1, i a_4, -i \varphi_3, i \varphi_2).
$$ 

Any solution $\chi$ to $\hD_{0,3} \chi = 0$ and such that $[\D_i,\chi]=0$, $i=1,2,3$ should be discounted since the corresponding perturbation 
$(a,\varphi)$ vanishes identically.  In fact, there are no viable infinitesimal gauge transformations of this form anyway. Indeed, recall that any commuting 
triple can be written as $\D_i=g\D_i^{(0)}g^{-1}$ for some complex gauge transformation $g$, where
\begin{equation}\label{mill}
\D_1^{(0)}=\partial_{\bar z}, \ \  \D_2^{(0)}=\partial_y,\ \ \mbox{and} \ \ \ \D_3^{(0)}=\begin{pmatrix}0 & z^\r\cr 0 & 0 \end{pmatrix},
\end{equation}
cf.\ \eqref{murky} and \eqref{urky}. The $\D_i$ also satisfy the moment map condition when $g$ is given by \eqref{turfew}.  We then calculate that 
the further perturbation \eqref{forf} vanishes to second order in $\varepsilon$, i.e., $[\D_i,\chi]=0$, $i=1,2,3$, precisely when 
 \begin{equation}\label{infico} 
\chi=g\begin{pmatrix} 0& f(z)\cr 0&0  \end{pmatrix}g^{-1}=\begin{pmatrix}0& e^v f(z)\cr 0 & 0 
\end{pmatrix}
\end{equation}
for some $f(z)$ which is holomorphic in $z$ and independent of $y$.  However, $e^v\sim 1/y$ as $y \searrow 0, \ z \neq 0 $;
however, we admit only gauge transformations which are bounded as $y \to 0$, so this choice is not admissible.

There is a small arithmetic difference in calculating the Type II indicial roots.  
By \eqref{sell}, the perturbations $a_i$ and $\varphi_i$ are first derivatives of $\chi$ or commutators of a component of $A$ or $\phi$ with $\chi$,
and these operations shift the exponent by $-1$. Thus for Type II perturbations, \eqref{weloxe} is replaced by  
\begin{equation}\label{elox}
\lambda=-\frac{3}{2}\pm \sqrt{\gamma+\frac{1}{4}}.
\end{equation}
The bound $\gamma>2$ now means that there are no Type II indicial roots in the closed interval $[-3,0]$.  For large $\r$, this result is again fairly sharp. 

\subsubsection{The Exceptions}\label{exceptions}
The calculation for Type II deformations is not yet complete, and to finish the analysis we must classify the Type II deformations which are
not of the form $\delta\D_i=[\D_i,\chi]$.   The indicial roots associated to these new deformations will be said to be of Type II$'$. 

As above, write $\D_i = g \D_i^{(0)} g^{-1}$ with $g$ as in \eqref{turfew}.  Up to a complex gauge transformation, any commutativity-preserving 
deformation can be assumed to modify only $\D_3^{(0)}$. We consider a general perturbation of this type: 
\begin{equation}
\label{usual}
\begin{aligned} 
\D_1& = g \frac{\partial}{\partial\bar z} g^{-1}, \qquad \D_2 = g\frac{\partial}{\partial y}g^{-1}\cr
\D_3&=g\left(\begin{pmatrix}0&z^\r\cr 0 & 0\end{pmatrix}+  \varepsilon 
\begin{pmatrix}\alpha & \beta\cr \gamma & -\alpha \end{pmatrix} + \calO(\varepsilon^2)\right) g^{-1}=
\begin{pmatrix} \varepsilon \alpha & (z^\r+ \varepsilon \beta)e^v \cr \varepsilon \gamma e^{-v}& - \varepsilon \alpha\end{pmatrix} + \calO(\varepsilon^2).
\end{aligned}
\end{equation}
We henceforth drop the terms of order $\varepsilon^2$ and higher. 
This perturbation preserves commutativity only if $\alpha,\beta,$ and $\gamma$ are holomorphic in $z$ and independent of $y$.   Next,
this perturbation is of type II, i.e., of the form $\delta\D_3= g [\D_3^{(0)},\chi] g^{-1}$ for some $\chi$, if and only if $\gamma=0$ and $\alpha$ and $\beta$ 
are divisible by $z^\r$. Thus a basis of commutativity-preserving perturbations of the $g \D_i^{(0)} g^{-1}$, modulo those of the form 
$\delta\D_3= g [\D_3^{(0)},\chi] g^{-1}$,  
is given by 
\begin{equation}\label{moro}
\begin{aligned} 
\alpha(z)&=z^t,~~0\leq t <\r\cr \beta(z)&=z^\sigma,~~0\leq \sigma <\r\cr \gamma(z)&=z^k,~~0\leq k. 
\end{aligned} 
\end{equation}

For $\alpha,\beta,\gamma$ as in \eqref{moro}, the perturbation \eqref{usual} satisfies neither the gauge condition $S=0$ nor the linearized moment map condition 
$\delta\mu=0$.  We therefore modify this linear perturbation further by a pure gauge term $[\D_i, \chi]$ for some complex-valued $\chi$.  The entire
linear perturbation is then 
\begin{equation}\label{welz}
[\D_1, \chi] , \ \ [\D_2, \chi], \ \ \ \mbox{and} \quad [\D_3, \chi]+  \begin{pmatrix} \alpha(z) & \beta(z)e^v \cr \gamma(z)e^{-v}& -\alpha(z)\end{pmatrix}.
\end{equation}
Since $\uphi^\dagger=\begin{pmatrix}0& 0\cr  \bar z^\r e^{v}& 0\end{pmatrix}$, then by \eqref{linSmu0}, the condition $- S + i \delta \mu = 0$ becomes
\begin{equation}\label{wolf}
\hD_{0,3} \chi=\left[\uphi^\dagger,   \begin{pmatrix} \alpha(z) & \beta(z)e^v \cr \gamma(z)e^{-v}& -\alpha(z)\end{pmatrix}\right] = 
\begin{pmatrix} -\bar z^\r e^{2v}\beta & 0 \cr 2\bar z^\r e^v\alpha & \bar z^\r e^{2v}\beta \end{pmatrix},
\end{equation}
where $\hD_{0,3}$ is the usual operator and $\alpha$, $\beta$ and $\gamma$ are any linear combinations of the terms in \eqref{moro}.

Homogeneous solutions of  $\hD_{0,3} \chi=0$ correspond to ordinary Type II deformations, which we already understand. Thus 
we focus on solutions of \eqref{wolf} satisfying the Nahm pole boundary condition away from $z=0$, modulo  solutions of the 
homogeneous equation. 

The right side of \eqref{wolf} does not have a fixed homogeneity in $\rho$, but it can be decomposed into homogeneous terms. 
First, if $\alpha=\beta=0$, this right hand side vanishes, so we may take $\chi=0$.  The corresponding perturbations 
$\delta\uphi=\begin{pmatrix}0 & 0\cr z^k e^{-v} & 0 \end{pmatrix}$ in \eqref{welz}  vanish at  $y=0$ since $e^{-v}\sim y$ as $y\searrow 0$, $z\neq 0$, 
hence are allowable. Since $e^{-v}$ is homogeneous in $\rho$ of degree $\r + 1$,  $\gamma(z)e^{-v}$ scales like $\rho^{k+\r+1}$. Thus 
the Type II$'$ indicial roots corresponding to these lower-triangular deformations of $\uphi$ are $\r+1,\r+2,\r+3,\dots$. 

Next suppose that $\alpha=0$, $\beta=z^\sigma$ (with $0\leq \sigma <\r$).  The right side of \eqref{wolf} is now homogeneous in $\rho$ of
degree $\sigma  - \r - 2$, so we let 
\begin{equation}\label{moft}
\chi = \rho^{\sigma - \r} \begin{pmatrix} f(\psi)e^{i(\sigma-\r)\theta} & 0 \cr 0 & -f(\psi)e^{i(\sigma -\r)\theta}\end{pmatrix}
\end{equation}
for some function $f(\psi)$. Inserting this into \eqref{wolf}, and writing $e^{2v} = \rho^{-2\r - 2} k(\psi)^2$, we see that 
$$
(M_S - (\r-\sigma)(\r -\sigma -1)) (e^{i(\sigma-\r)\theta} f(\psi) H ) =  e^{i(\sigma-\r)\theta} k(\psi)^2 H, \qquad  H = \begin{pmatrix}  1 & 0 \\ 0 & -1 \end{pmatrix}.
$$
(We could of course reduce even further to an ODE in $\psi$.) Since $\rho^{2\r + 2} e^{2v} = k(\psi)^2 \sim 1/\cos^2 \psi = 1/s^2$, we can apply Proposition~\ref{MSprop}
to conclude that there exists a solution $f \in \calC^\infty(S^2_+)$.   In order that the Nahm pole singularity persist in this perturbation, the leading terms in the
two summands of $\delta \D_3$ must cancel at $\psi = \pi/2$.  Using that
$$
\D_3 = \rho^\r \sin^\r \psi \, e^{i\r \theta} e^v X,\ \ \chi = \rho^{\sigma - \r} e^{i(\sigma - \r)\theta}  f(\psi) X, \ \mbox{where}\ X = \begin{pmatrix} 0 & 1 \\ 0 & 0 \end{pmatrix},
$$
this condition becomes 
$$
\rho^\sigma \sin^\sigma \psi \, e^{i\sigma \theta} e^v (-2f ) X + \rho^\sigma \sin^\sigma \psi \, e^{i\sigma \theta} e^v X \sim 0,\ \ \mbox{as}\ \psi \to \pi/2,
$$
or finally, $f(\pi/2) = 1/2$. 

It remains to show that $(\r-\sigma)(\r-\sigma-1)$ is not in the spectrum of $M_S$ and that $f = 1/2$ at $s=0$. For the first part, note that on diagonal
matrices, $D_\theta$ reduces to $\del_\theta$, so the action of $M_S$ on $e^{i(s-\r)\theta} f H$ reduces to
$$
(- \frac{\del^2\,}{\del s^2} + \frac{\sin s}{\cos s} \frac{\del\,}{\del s} +  (\sigma-\r)^2 + N_S) (fH) = -k^2 H. 
$$
Now $N_S > 0$, and the other terms on the left consitute the ordinary scalar Laplacian acting on functions $f(\psi) e^{i(s-\r)\theta}$ on $S^2_+$. The 
smallest eigenvalue of the scalar Laplacian on such functions, even on the entire sphere, is $(\r-\sigma)(\r-\sigma -1)$, so the smallest eigenvalue of $M_S$ must
be even larger.  Finally, since $N_S \sim 2/s^2$ and $k^2 \sim 1/s^2$ as $s \to 0$, we may compute formally, using the smoothness of $f$ at $s=0$, to see that 
$f = 1/2$ there. 

Recalling that $(a,\varphi)$ is homogeneous of one degree lower than $\chi$, we have obtained new Type II$'$ linearized solutions of the KW equations 
not seen in section \ref{generic}, with negative indicial roots $\sigma-\r-1$, $0\leq \sigma <\r$. 

For the final case, set $\beta=0$ and $\alpha=z^t$, $0\leq t<\r$.  Proceeding as before,  since $e^v \bar{z}^{\r} \alpha = \rho^{t-1} e^{i(t-\r)\theta} k(\psi)$, 
we search for $\chi$ of the form 
 \begin{equation}\label{loft}
\chi=\rho^{t+1} e^{i(t-\r)\theta} \begin{pmatrix}0 & 0\cr  h(\psi) & 0\end{pmatrix},
\end{equation}
so that \eqref{wolf} becomes 
$$
(M_S - t (t+1)) e^{i(t-\r)\theta}  (h(\psi) Y)  = 2 e^{i (t-\r)\theta} k(\psi) Y, \qquad  Y = \begin{pmatrix} 0 & 0 \\ 1 & 0 \end{pmatrix}.
$$
Since $[H, Y] = -2Y$, $D_\theta$ reduces to $\del_\theta + i r\del_r v$, so equivalently 
\begin{equation}\label{wolf2}
(- \frac{\del^2\,}{\del s^2} + \frac{\sin s}{\cos s} \frac{\del\,}{\del s} +  (t-\r + r\del_r v)^2 -t (t+1) + N_S) (h(s) Y) = k(s) Y.
\end{equation}

As before, we need an eigenvalue estimate to show that the operator on the left is invertible, but this is no longer quite as direct.
We first compute that 
$$
N_S ( Y) = (|\mu|^2 + |\nu|^2) Y, 
$$
where
\begin{equation*}
\begin{aligned}
\mu & = - (\r + 1) \frac{ (1 + \cos \psi)^{\r + 1} + (1 - \cos \psi)^{\r + 1}}{ (1 + \cos \psi)^{\r + 1} - (1 - \cos \psi)^{\r + 1}} \\
\nu & = 2 (\r + 1)  \frac{(\sin \psi)^\r e^{i\r \theta}}{ (1 + \cos \psi)^{\r + 1} - (1 - \cos \psi)^{\r + 1}}.
\end{aligned}
\end{equation*}
whence, reverting back to $s = \pi/2 - \psi$,
$$
M_S - t(t+1) \geq -\frac{\del^2\,}{\del s^2 } + \frac{ \sin s}{\cos s}\frac{\del \,}{\del s} + |\mu|^2 - t(t+1) = -\frac{1}{\cos s} \frac{\del\,}{\del s} 
\left( \cos s \frac{\del\, }{\del s}\right) + |\mu|^2 - t(t+1).
$$
The nonnegative terms $- (\cos s)^{-2} D_\theta^2$ and $|\nu|^2$ have been dropped. 
We conclude invertibility of the operator in \eqref{wolf2} by noting that $|\mu|^2 \geq (\r + 1)^2 > t(t+1)$,
which follows by the simple observation that $|\mu| = (\r + 1) (G+F)/(G-F)$, where $0 \leq F < G$

We may now invoke Proposition~\ref{MSprop} again to solve for $h$. It remains to study the asymptotics of this solution as $s \to 0$.  Recall that the indicial 
roots of $M_S$ (and $M_S - t(t+1)$) equal $-1$ and $2$.  Since $k \sim 1/s$, this same Proposition implies that $h \in s \calC^\infty(S^2_+)$.  This is
already enough to ensure that $\delta \D_3 = \calO(1)$ as $s \to 0$, so this perturbation does not compete with the Nahm pole singularity. However, we can
do slightly better: inserting $h \sim a s + \calO(s^2)$ into the defining equation and recalling that $N_S \sim 2/s^2$, gives that $h(0) = -1$.  Finally,
$$
\delta \D_3 = [ z^\r e^v X, \rho^{t+1} e^{i(t-\r)\theta} h Y] + z^t H \sim \rho^t e^{it\theta} ( s^{-1}h(0)  + 1) + \calO(1), \ \ \mbox{as}\ s \to 0,
$$
so in fact, $\delta \D_3 = \calO(s)$.  

This produces the final set of Type II$'$ indicial roots $\{0, 1, \ldots, \r - 1\}$. 

\bigskip

We summarize all of this in the following
\begin{prop}\label{fullindrts}
The set of indicial roots of the linearized KW equations are as follows
$$
\begin{aligned}
&\mbox{Roots of Type I}: ~{\mathrm{a~ subset ~of}}~  \{ -\frac12 \pm \sqrt{\gamma + \frac14}:  \gamma \in \mathrm{spec}\,(M_S) \}; \\
& \mbox{Roots of generic Type II}:   \{ -\frac32 \pm \sqrt{\gamma + \frac14}:  \gamma \in \mathrm{spec}\,(M_S) \}; \\[0.5em]
& \mbox{Roots of Type II$\,'$}:  \{-\r-1, -\r, \ldots, -2\} \cup \{0, 1, 2, \ldots, \r - 1\} \cup \{ \r+1, \r+2, \ldots \}.
\end{aligned}
$$
Listed for Type I are the indicial
roots for $\LKW^\dag \LKW$; some may  not be  indicial roots for $\LKW$ alone.

Amongst these, the roots of Type I correspond to perturbations with $a_1$ or $\varphi_4$ nonzero, 
while those of generic Type II and Type II$\,'$ correspond to perturbations with  $a_1 = \varphi_4 = 0$. 

The eigenvalues $\gamma\in \mathrm{spec}\,M_S$ are all greater than 2, so that there are no Type I roots in the interval $[-2,1]$ and no
generic Type II roots in the interval $[-3,0]$.
\end{prop}

\medskip

The indicial root $0$ of Type II$'$ is particularly interesting because of the symmetries of the corresponding fluctuation. This fluctuation 
is bounded at $y=z=0$, hence is one power of $\rho$ less singular than the model solution at a knot.  It satisfies the reduced KW equations 
to first order, and in this framework, which involves the commuting operators $\D_1, \D_2, \D_3$, there are independent symmetries 
that rotate the complex $z$-plane or multiply $\uphi$ by a complex number of modulus $1$. The rotation symmetry of the full 
KW equations is the diagonal combination of these two symmetries;  the full KW equations do not have these two separate symmetries. 
The diagonal symmetry acts naturally on the $(1,0)$-form $\uphi\,\d z$ 
in the complex $z$-plane.  Setting $\alpha$ to be a constant, say $w$, the perturbation of the $(1,0)$-form is
\begin{equation}\label{loro} 
\delta(\uphi\,\d z)=\begin{pmatrix} w & 0 \cr 0 & -w \end{pmatrix}\d z, 
\end{equation}
which is not invariant under rotations in  $z$, but transforms with angular momentum $1$.  The deformation space of the reduced
KW equations has a natural complex structure, and using this it is natural to think of $w$ as a complex parameter.  But there is no 
natural complex structure for the solution space of the full KW equations, so one should then think of $w$ as a pair of real 
parameters -- its real and imaginary parts -- that transform under rotations as the ``vector'' of
$SO(2)$, corresponding to  angular momentum $\pm 1$.

\subsection{The Nahm pole boundary condition} 
Following our comments at the end of section 2, we may now give the `quantitative' description of the generalized Nahm pole boundary 
condition for the KW equations on a Riemannian manifold with boundary $(M,g)$ with knot $K \subset W = \del M$. This involves 
two pieces of data: the first is an injective bundle map $\phi_\rho: TW_{W \setminus K} \to \ad (E)$ associated to a principal representation 
$\rho: \mathfrak{s}\mathfrak{u}(2) \to \mathfrak{g}$, and the second is the field $(A^K, \phi^K)$ defined along $K$ where, in 
Fermi coordinates around this knot, $A^K$ and $\phi^K$ equal the model solution given in \eqref{more}.  Due to the rotational 
symmetry in $z$ and the (built-in) independence of $x_1$, there is no ambiguity in the leading order term in the expansion for 
$(A^K, \phi^K)$ as $\rho \to 0$.  As already indicated, solutions of the full KW equations 
are fields $(A,\phi)$ which differ from these leading order models by lower order perturbations.

We use both coordinate systems $(x, y)$ and $(\rho, s, \theta, t)$.  As in \cite{MW}, away from $K$ we search for solutions 
$(A,\phi) = (A^K, \phi^K) + (a,\varphi)$ with $|a|, |\varphi| = \calO(\rho^{-1 + \varepsilon} s^{-\varepsilon})$ and 
$|\varphi| = \calO( \rho^{-1+\varepsilon} s^{-1 + \varepsilon})$. We have also already checked the compatibility of these 
conditions in the region where $\rho \to 0$ and $s \to 0$.  The rates of vanishing or blowup
here are chosen so that these fluctuations do not interfere with the leading order terms. These will be formalized even further later
on by requiring that the fluctuations lie in certain weighted H\"older spaces. 

The analytic work later in this paper justifies the fact that solutions of the KW equations satisfying these boundary conditions 
admit complete asymptotic expansions at $W \setminus K$ and at $K$ -- a regularity condition called polyhomogeneity.
That analysis provides the step intermediating between the formal rates of decay of linearized solutions as described by our
anaysis of solutions of the indicial equation and the actual rates of decay (and higher regularity) for solutions of the linearized
and nonlinear KW equations.  In particular, as indicated by the results in section \ref{indicial},  there are no indicial roots in the semi-open 
interval $[-1,0)$, and no indicial roots of Type I, i.e., with $a_1$ or $\varphi_4$ non-zero, in the larger interval $[-2,1]$. This will be used to
show that the perturbation $(a,\varphi)$ about $(A^K, \phi^K)$ is bounded as $\rho \to 0$ and in accordance with the behavior
of fluctuations away from $K$, decays as $s \to 0$. 

In the next sections we describe two applications of this regularity theory for solutions of the KW equations which satisfy
the generalized Nahm pole boundary conditions.  The first is a Weitzenb\"ock formula in the presence of knots which generalizes 
the one in \cite{MW} when $K$ is absent. The second is an index formula for $\LKW$. 

\section{The Weitzenb\"ock Formula With Knots}\label{weitzknots}

Let us now return to the model setting $X = \RR^4_+$ with $K = \RR$ and describe the Weitzenb\"ock  formula which generalizes the
one in \cite{MW} to the present setting.  

We begin by recalling the basic Weitzenb\"ock  formula, eqn.\ (2.8) in \cite{MW}.  First write the KW equations as 
\begin{equation}\label{elboxoc} \V_{ij}=\V^0=0, \end{equation}
where
\begin{equation}\label{telbo}\V_{ij}=F_{ij}-[\phi_i,\phi_j]+\epsilon_{ij}{}^{kl}D_k\phi_l,~~~~\V^0=D_i\phi^i. \end{equation}
On a closed manifold, one calculates that
\begin{equation}\label{zoffbo} -\int_M\d^4x \sqrt g\Tr\left(\frac{1}{2}\V_{ij}\V^{ij}+(\V^0)^2\right)=I,
\end{equation}
where 
\begin{equation}\label{baction}
I=-\int_M\d^4x \sqrt g\Tr\left(\frac{1}{2}F_{ij}F^{ij}+D_i\phi_j D^i\phi^j+R_{ij}\phi^i\phi^j+\frac{1}{2}[\phi_i,\phi_j][\phi^i,\phi^j]\right); \end{equation}
in these formul\ae, $R_{ij}$ is the Ricci tensor and summation over repeated indices is understood. This identity can be written in a coordinate-invariant
way and requires no assumptions about the geometry. Together with its generalization to $t\not=\pm 1$, it leads to powerful vanishing theorems for the 
KW equations on a four-manifold $X$ without boundary. 

If $X$ has boundary with Nahm pole boundary conditions, \eqref{baction} must be emended to 
\begin{equation}\label{zoffbo2} -\int_M\d^4x \sqrt g\Tr\left(\frac{1}{2}\V_{ij}\V^{ij}+(\V^0)^2\right)=I + \int_{\partial M}\d^3x\, 
\epsilon^{abc}\Tr\left(\frac{1}{3}\phi_a[\phi_b,\phi_c]-\phi_a F_{bc}\right). \end{equation} 
(We write $i,j,k=1,\dots,4$ for indices tangent to $X$ and $a,b,c=1,\dots,3$ for indices tangent to $\partial X$.)  However, the Nahm singularity makes
these boundary contributions infinite, so this Weitzenb\"ock  formula is not very useful.  

To counter this, we introduced a second Weitzenb\"ock  formula adapted to the Nahm pole condition \cite[Eqn.\ (2.21)]{MW}. To find this formula, one defines an
alternate expression $I'$ for the bulk integral which is a sum of squares of certain quantities that vanish in the basic Nahm pole solution on a half-space, and
whose vanishing characterizes this basic solution. The tricky point which made the analysis in \cite{MW} possible is the not-so-obvious identity \cite[Eqn.\ (2.21)]{MW}, 
which says that $I'$ equals the sum of squares of the KW equations plus an exact term.   

In the presence of a knot \cite[Eqn.\ (2.21)]{MW} must be refined for the analogous reason:  the boundary terms there diverge because of the extra singularities 
near the  knot.  To handle this, we replace $I'$ by another functional $I''$ which is again a sum of squares of certain quantities which are chosen to vanish for 
the model solution with a knot.  

To determine these quantities, we return to the reduced KW system on $\RR^3_+$. Let us set
\begin{equation}\label{zub}
\X=[\D_1,\D_2],~~ \Y=[\D_2,\D_3],~~\Z=[\D_3,\D_1]. 
\end{equation}
Solutions of the reduced system are characterized by 
\begin{equation}
\label{ub}
\begin{aligned}\X=\Y=\Z=\mu &= 0  \cr
F_{1a}=D_1\phi_a=D_a\phi_4=[\phi_a,\phi_4]&=0,~~1\leq a \leq 4.
\end{aligned}
\end{equation}
The conditions in the second line say that, up to a gauge transformation, the solution is independent of $x^1$, and moreover that $\phi_4$ is covariantly
constant and generates a symmetry of the solution.  Supplementing \eqref{ub} by the condition that $\phi_4$ vanishes somewhere on a finite or infinite boundary  of
$\R^4_+$ (in our applications it will vanish on all finite and infinite boundaries), we obtain that $\phi_4 \equiv 0$ so we have a solution of the reduced KW equations.                            
The characterization (\ref{ub}) of the reduced KW equations suggests that we define 
\begin{equation}\label{tub}
I''=-\int_{\R^4_+}\d^4x\,\Tr\left(\X\bar\X+\Y\bar\Y+\Z\bar \Z+\mu^2+\sum_{a=1}^4\left(F_{1a}^2+D_1\phi_a^2+[\phi_a,\phi_4]^2\right) +
\sum_{b=2}^4D_b\phi_4^2 \right).
\end{equation} 
Assuming $\phi_4$ is known to vanish somewhere on $\partial \R^4_+$, the condition $I''=0$ characterizes pairs $A,\phi$ that are gauge-equivalent
to a solution of the reduced KW equations.  

This definition of $I''$ is useful because of the following analogue of \cite [Eqn.\, (2.21)]{MW}: 
\begin{equation}\label{welub}
-\int_{\R^4_+}\d^4x \left(\frac{1}{2}\V_{ij}\V^{ij}+(\V^0)^2\right)=I''+\Omega, 
\end{equation} 
where $\Omega =\Omega_1+\Omega_2+\Omega_3$, with
\begin{equation}
\label{telub}
\begin{split}
\Omega_1=-2 \int_{\R^4_+}\d^4x&\,\left( \sum_{a=1}^3\left(\partial_a\Tr \phi_a D_4\phi_4-
\partial_4\Tr\,\phi_aD_a\phi_4\right)  +\sum_{a=2}^3\left(\partial_1\Tr\,\phi_1D_a\phi_a -\partial_a\Tr\,\phi_1D_1\phi_a
\right)\right) \\
\Omega_2=-2 \int_{\R^4_+}\d^4x&\, \biggl(-\partial_1 \Tr\, \phi_2[\phi_3,\phi_4]+\partial_2\Tr\,\phi_3[\phi_4,\phi_1]-\partial_3\Tr\,\phi_4[\phi_1,\phi_2]  \biggr)\cr
\Omega_3=-2 \int_{\R^4_+}\d^4x& \,\biggl( \partial_1\Tr\,\left(\phi_2F_{34}+\phi_3F_{42}+
\phi_4 F_{23}\right)+\partial_2\Tr\,\left(-\phi_3 F_{41}-\phi_4F_{13}\right)\biggr. \\ 
&  - \biggl.\partial_3\Tr\,\left(\phi_2F_{41}-\phi_4F_{12}\right)+\partial_4\Tr\,\left(-\phi_3 F_{12}+\phi_2F_{13} \right) \biggr) .
\end{split}
\end{equation}

If a solution of the KW equations on $\RR^4_+$ has sufficient regularity at $y=0$, and in particular along the knot $\RR$, 
and sufficient decay at infinity, then  integration by parts in $\Omega_1,\Omega_2,\Omega_3$  shows that these
boundary contributions vanish. We investigate this here using the behavior formally predicted by the indicial root computations above,
and analogous computations near infinity. The justification is a focus in the later part of this paper. 

There are three separate regions where the boundary behavior  needs to be examined: (1) as $y\to 0$ away from the knot; 
(2) along the knot, with $x_2, x_3, y \to 0$ simultaneously; and (3) as $(x_1, x_2, x_3, y) \to \infty$.  (As always, we use
$y$ and $x_4$ interchangeably.) 

Case (1) was already treated in \cite{MW}. That paper contains a simpler version of the Weitzenb\"ock  formula, and the boundary 
terms arising in that formula are simpler than those in \eqref{telub}.  In this region, by virtue of the polyhomogeneity (tangential
regularity and asymptotic expansions in the normal direction) of solutions, we may disregard the terms that involve differentiations 
with respect to $x_1, x_2, x_3$ since these do not alter the rates of blowup or decay. The terms in \eqref{telub} which involve 
derivatives with respect to $y$ coincide with the analogous terms in eqn. (2.21) of \cite{MW}. We may then follow the arguments in 
that paper to show that these boundary terms vanish. 

For case (2), the boundary terms along the knot, denote by $H_\varepsilon$ the hemisphere of radius $\varepsilon$ 
around $z=y=0$ at fixed $x_1$, which we consider as lying in the boundary of $\R^3_+$. The actual boundary term
in $\RR^4_+$ involves integrating over $\RR \times H_\varepsilon$, where $\RR$ is the $x_1$ direction, and so long
as we are integrating along a finite interval in $\RR$, it suffices by Fubini to prove vanishing of the boundary term
along the hemisphere as $\varepsilon \to 0$.

We first verify that these boundary terms actually vanish for the model solution itself.  This is not completely trivial 
since the area of $H$ is of order $\varepsilon^2$, whereas counting powers suggests that the integrand in this boundary
term is of order $1/\varepsilon^3$ (indeed, $\phi$, $D\phi$, and $F$ are of orders $1/\varepsilon$, $1/\varepsilon^2$, and 
$1/\varepsilon^2$, respectively, so $\phi^3$, $\phi D\phi$, and $\phi F$ are all of order $1/\varepsilon^3$).  However, fortunately
$\phi_4=F_{1i}=D_1\phi_i=0$ in the model solution, which yields the desired vanishing ($\Omega_1$ causes no trouble because 
$\phi_4=D_1\phi_a=0$, $\Omega_2$ because $\phi_4=0$, and $\Omega_3$ because $F_{1a}=0$; surface terms involving
derivatives with respect to $x^1$ are not relevant here). 

For the next part, write $\Psi$ schematically for the model solution $(A, \Phi)$ and $\delta \Psi$ for the fluctuation 
term $(a,\varphi)$. The polyhomogeneity of $\delta \Psi$ is proved in section 8, and the precise decay rates of these terms are those
predicted by the indicial roots calculated earlier.  We then insert $\Psi + \delta \Psi$ into the boundary terms.

Consider first the contributions from the parts that are linear in $\delta\Psi$. For the moment, suppose
this involves the integral over $H$ of a term of the form $\Psi^2\delta\Psi$.  If $\delta\Psi \sim \varepsilon^\lambda$, 
then $\Psi^2\delta\Psi \sim \varepsilon^{-2+\lambda}$, so $\int_H \Psi^2\delta\Psi \to 0$ provided $\lambda>0$
The only possible problem is the exceptional Type II$'$ mode for which $\lambda=0$. 
However, this mode has $\delta \phi_4= \delta A_1=0$; in particular the vanishing of $\delta\phi_4$ ensures that 
$\int_H \Psi^2\delta\Psi \to 0$.

We must also consider terms $D\Psi\delta\Psi$, $\Psi D\delta\Psi$, and $F\delta\Psi$.  These scale in the same way as 
$\Psi^2\delta\Psi$ except that terms $D_1\Psi\delta\Psi$ or $\Psi D_1\delta\Psi$ have an extra factor of $\varepsilon$,
since (using tangential regularity and that $A_1=0$ in the model solution) a covariant derivative along the knot does not 
increase the singularity.  Because we only consider perturbations with `decay' rate $\lambda > -1$, these terms are
harmless too.  As before, counting powers shows that the modes with $\lambda>0$ do not contribute, and the 
exceptional mode with $\lambda=0$ is also not an issue because $\delta\phi_4=\delta A_1=0$ eliminates all boundary terms 
near the knot except for those which drop out because $D_1$ does not increase the singularity.

We have only considered contributions linear in $\delta\Psi$. The terms quadratic or cubic in $\delta\Psi$ are only 
less singular. Hence all boundary terms in case (2) vanish. 

For case (3), we proceed very much as in \cite[Section 2.6]{MW}. The idea is simply that by hypothesis, $\delta \Psi$ tends
to infinity as $R = |(x_1, x_2, x_3, x_4)| \to \infty$, and we must verify that this rate is such that if we integrate
by parts on some large region, say $R \leq R_0$,  then the boundary terms tend to zero as $R_0 \to \infty$.  To 
prove this decay, we employ the reasoning explained in the beginning of section 4.1: consider a solution 
$\Psi + \delta \Psi$ of the KW equations in the region where $R$ is large, so that $\delta\Psi$ is small, and 
expand the action $I$ into terms of order $\delta\Psi$, $|\delta \Psi|^2$, etc.  The linear term vanishes because $\Psi$ 
is already a solution, and the stationarity of the perturbed solution implies that $\calL^\dag \calL \delta \Psi$
must vanish modulo terms of higher order in this asymptotic regime. This means that the decay rate of solutions 
of the linear equation $\calL^\dag \calL (a,\varphi) = 0$ determine the decay rate of the nonlinear fluctuation $\delta \Psi$. 
To analyze this we use the linear Weitzenb\"ock  formula $\calL^\dag \calL = \hD$.  

As described in section 4.1, since $\calL$ is the linearization at the model knot solution, this operator does not couple 
$a_1$ or $\varphi_4$ with the other terms, and so reduces to 
$$
\hD_{0,3} a_1 = \left(-\sum_{i=1}^4 D_i^2 + \sum_{j=1}^3 [\phi_j, [\phi_j, \cdot]]\right) a_1,
$$
and similarly for $\varphi_4$.    Assume this expression vanishes. Decomposing $a_1$ into a sum of eigenfunctions 
for the operator $M_S$ on $S^2_+$, cf.\ \eqref{welox}, each eigencomponent $a_{1,\gamma}$ satisfies the equation 
$$
\left(- \frac{\del^2 \,}{\del \rho^2} - \frac{\del^2\,}{\del x_1^2} - \frac{2}{\rho} \frac{\del\,}{\del \rho} + 
\frac{\gamma}{\rho^2}\right) a_{1,\gamma} = 0,
$$
where $\gamma$ is the eigenvalue of $M_S$.  This may be transformed further, using polar coordinates
$\rho = R\sin \omega$, $x_1 = R \cos \omega$, to the equation
\begin{equation}
\left(-\frac{\del^2\,}{\del R^2} - \frac{3}{R}\frac{\del\,}{\del_R}  + \frac{1}{R^2}( -\frac{\del^2}{\del \omega^2} + 
\cot \omega \frac{\del\,}{\del \omega} + \gamma) \right)a_{1,\gamma} = 0.
\label{eigeneqn}
\end{equation}
The operator involving $\omega$ is self-adjoint with respect to the measure $(\sin \omega)^{-1} d\omega$, and
is bounded below by $\gamma$, so the decay rate of solutions is better than $-1 - \sqrt{1 + \gamma} < -2$.
(In fact, $\gamma > 2$ and $-\del_\omega^2 + \cot \omega \del_\omega \geq c > 0$ on functions which vanish at 
$\omega = 0, \pi$, so this decay rate can be improved.) We conclude that $a_1$ and $\varphi_4$ decay like $R^{-2-\varepsilon}$ 
for some $\varepsilon > 0$.

To understand decay rates of the other components, we must analyze the full operator $\hD = \LKW^\dag \LKW$,
which appears in section 4.1 as the Euler-Lagrange operator for the action $I_{2,0} + I_{2,1}$ in \eqref{ztwo}.  
We now must study this operator acting on the remaining fields $(a', \varphi') = (a_2, a_3, a_4, \varphi_1, \varphi_2, \varphi_3)$
We deduce from \eqref{ztwo} that 
$$
\hD(a',\varphi') = - \left(\sum_{i=1}^4 D_i^2 + \sum_{j=1}^3 [\phi_i, [\phi_i, \cdot] ]\right) (a',\varphi') + P (a',\varphi'),
$$
where 
$$
P(a',\varphi') = - 2 \sum_i \left( ( [F_{ij}, a_i] + [D_j \phi_i, \varphi_i],  [D_i \phi_j, a_i] + [ [\phi_i, \phi_j], \varphi_i] \right).
$$
Note that each coefficient of $P$ decays like $1/\rho^2$ as $\rho \to \infty$. Indeed, writing $P = \rho^{-2} P_S$, and
referring back to the notation \eqref{welox}, we have
\begin{equation}
\hD = -\frac{\partial^2}{\partial x_1^2} - \frac{\partial^2}{\partial\rho^2} - \frac{2}{\rho}\frac{\partial}{\partial\rho}
- \frac{1}{\rho^2}(\Delta_{A} + N_S + P_S),
\label{NSPS}
\end{equation}
where $\Delta_A$ is the connection Laplacian on $S^2_+$ (this restriction is valid since $A$ only has a $d\theta$ component),
$N_S$ is defined in \eqref{omolt}, and we write $N_S + P_S = Q_S$ for convenience.  Note that $Q_S \sim C Q_S' \cos \psi^{-2} = 
C Q_S' s^{-2}$ as $s \to 0$. This is enough to ensure that $\Delta_A + Q_S$ has discrete spectrum, see Proposition \ref{MSprop}. 
Since $\hD$ comes from a nonnegative quadratic form, $\Delta_A + Q_S$ is also a positive operator.  

We can now apply the same line of reasoning as before.  First write $\hD$ in cylindrical coordinates $(\rho, x_1, \Theta)$,
$\Theta \in S^2_+$.  Let $\gamma$ be any eigenvalue of $\Delta_A + Q_S$. Then the Fourier coefficients $(a', \varphi')_\gamma$ 
of this decomposition satisfy \eqref{eigeneqn}, though now $\gamma$ is an eigenvalue of $\Delta_A + N_S + P_S$
instead of just $\Delta_A + N_S$.  As noted, $\gamma > 0$, so we conclude, just as before, that these
coefficients decay like $R^{-2-\varepsilon}$ for some $\varepsilon > 0$. 

Altogether, these decay conditions are sufficient to ensure that the boundary term along $S_R = \{|(x_1, x_2, x_3, x_4)| = R\}$ 
tends to $0$ as $R \to \infty$. Indeed, the boundary terms are either cubic, or else the product of a component of the field 
$\phi$ or the solution $\varphi$ with a curvature two-form. The analysis above shows that such product decays like $R^{-3-\epsilon}$,
and in most cases faster, whereas these spheres have volume $R^3$, so the  boundary integrals vanishs as $R \to \infty$.  
The same argument as earlier shows that the extra potential contribution at the intersection of $S_R$ with the knot vanishes.

\section{Index}\label{index}
In this brief section we rely on the results of  Sections \ref{analKW} and \ref{param}, where it is proved that the linearized KW operator $\LKW$ 
with generalized Nahm pole boundary conditions is Fredholm between certain adapted function spaces.  We invoke this 
here to explain that the computation of the index of $\LKW$ in the knot-free case in \cite{MW} can be carried over to the 
present setting.
\begin{prop}
Let $(M,g)$ be a compact four-manifold with boundary, with a knot $K \subset \del M$.  Let $(A,\phi)$ be a solution to the 
KW equations satisfying generalized Nahm pole boundary conditions at $\del M$ with a knot singularity at $K$.  Denote by
$\LKW$ the linearization $\LKW$ of these equations at this pair of fields. Then 
\begin{equation}
\LKW: L^2(M, (\Lambda^1 \oplus \Lambda^3)\otimes \mathrm{ad}\, \mathfrak g) \longrightarrow L^2(M, 
(\Lambda^0 \oplus \Lambda^2 \oplus \Lambda^4)\otimes \mathrm{ad}\, \mathfrak g), 
\label{LKWfull}
\end{equation}
is an unbounded Fredholm operator.  If the metric $g$ is cylindrical near $\del M$, then 
$$
\mathrm{index}\,(\LKW) = - 3 \chi(M). 
$$
\end{prop}
We sketch the proof here since it is essentially the same as in \cite{MW}.  We prove in Section 8 below that this mapping 
is Fredholm, so we turn immediately to the calculation of its index when $g$ is cylindrical near the boundary.  
The key observation in \cite{MW} is that in the very special case that $M$ is an exact metric cylinder $W \times [0,1]$,
then $\LKW$ enjoys a symmetry which implies that its index vanishes. This symmetry holds even if we impose different 
Nahm pole boundary conditions at each end, including when either boundary component of this cylinder contains a knot 
and we impose the generalized Nahm pole boundary conditions there, or, in the other extreme, when one of the boundary 
components contains no knot and we impose either the classical relative or absolute boundary conditions at that boundary. 
We then use a standard excision argument to reduce the calculation of the index to the case where, for example,
the boundary condition is one of these classical ones, in which case the index calculation follows from standard Hodge
theory on manifolds with boundary.    If the metric is not cylindrical near the boundary, we should expect an index
formula with an extra local boundary correction term. 

In slightly more detail, if $g$ is cylindrical near $\del M$ and we consider fields with a knot singularity along $K \subset W = \del M$,
then we decompose $M$ into $M_1 \cup M_2$, where $\del M_1 = W$ and the restriction of $g$ to $M_1$ is still cylindrical near the 
boundary, and where $M_2 = W \times [0,1]$ is a metric cylinder.  Now impose relative boundary conditions at $\del M_1$,
and also at one end of $M_2$, and the original generalized Nahm boundary conditions at $W \times \{1\}$.  The excision 
principle shows that
$$
\mathrm{index}\,(\LKW, M) = \mathrm{index}\, (\LKW, M_1) + \mathrm{index}\,(\LKW, M_2),
$$
(where $(\LKW, M)$ of course means that we consider $\LKW$ on all of $M$, etc.).  The symmetry of $\LKW$ on $M_2$ alluded
to earlier shows that this last term vanishes, while the former equals $-3 \chi(M)$, see \cite [Section 4.2]{MW}.

It remains then to explain this symmetry, which was called pseudo skew-adjointness in \cite{MW}.  Its existence relies on
the distinguished direction along the generator of the cylinder, which we label by the variable $y$. Then with $a = 1, 2, 3$, 
the tangent fields $(a, \varphi)$ have components $(a_a, a_y, \varphi_a, \varphi_y)$. Define the endomorphism
$$
N \begin{pmatrix} a_a \\ \varphi_y \end{pmatrix} = \begin{pmatrix} \varphi_a \\  a_y \end{pmatrix}, \quad 
N \begin{pmatrix} \varphi_a \\ a_y \end{pmatrix} = - \begin{pmatrix} a_a, \\\varphi_y \end{pmatrix}. 
$$
Clearly 
$$
N^2 = -1, \ \ \mbox{and}\ \ \ N^\dag = -N,
$$
and it is explained in \cite{MW} that
$$
\LKW^\dag = -N \LKW N^{-1}. 
$$
Here $\LKW^\dag$ is the $L^2$ adjoint defined relative to the standard inner product 
$$
- \Tr \sum_{i=1}^4 (a_i^2 + \varphi_i^2)
$$
and volume form, both induced by $g$. More precisely still, this is the adjoint of the (unbounded) Fredholm map \eqref{LKWfull}.  
Since $N$ is an isometry on $L^2(M, \Lambda^*\otimes \mathrm{ad}\, \mathfrak g)$, it identifies the kernel and cokernel of
\eqref{LKWfull}, and hence shows that the index on this cylinder vanishes. 

In this entire discussion, the important point is that letting $\LKW$ act on $L^2$ already respects the Nahm pole boundary 
conditions, with or without knots.   As we have explained and the indicial root calculations make clear, the fluctuations
$(a, \varphi)$ must be `lower order' than the principal terms of the approximate solution $(A^K, \phi^K)$.  However, $A^K$ is 
bounded on $W$ away from $K$ while $\phi^K$ is not, so it seems that we should actually require $a$ to vanish at $W \setminus K$
while only imposing that $\varphi = \calO(y^{-1+\epsilon})$ there; on the other hand, both $a$ and $\varphi$ may be allowed
to be $\calO( \varrho^{-1+\epsilon})$ along $K$.  This may seem problematic because $N$ interchanges the components of these
two fields.  The resolution to this is the regularity theorem proved in \cite{MW}, that if $\LKW(a,\varphi) = 0$ even locally
in some neighborhood $\calU$ in $M$ around a point $p \in W$, then $(a,\varphi)$ have polyhomogeneous expansions which
in particular guarantees that $a$ decays as $y \to 0$.  In other words, the nullspace of $\LKW$ on $L^2$ is precisely the
tangent space to the relevant moduli space of solutions to ${\KW}(A,\phi) = 0$.  

\section{Compact knots on the boundary of $\RR^4_+$}\label{cpctK}
A particular setting of interest for this theory is when $K$ is a closed knot lying in the boundary of $\RR^4_+$.  
We discuss here briefly the analytic issues in this setting. 

We consider fields $\Psi = (A,\phi)$ which satisfy the Nahm pole boundary conditions along $y=0$,
with the knot singularity along $K$. We may compare these to the model Nahm solution $\Psi_0 = (0, \phi_\varrho)$ 
in the far field, and choose these to satisfy $|\Psi -\Psi_0| \leq C R^{-2-\varepsilon}$ for $R \geq R_0$. 
This condition implies that $\Psi - \Psi_0 \in L^2(\RR^4_+)$. As discussed in \cite{MW}, this decay is reasonable 
if $\Psi$ has only a Nahm pole singularity but no knot.   We claim that the same decay rate holds for 
fields with a singularity along a compact knot. The proof may be transferred verbatim from 
\cite[section 2.6]{MW}, since that discussion is local in such exterior regions. 
It is also the case that the index of $\LKW$ vanishes in this setting.  The proof of section \ref{index} may
be applied directly to this case.

\section{Analysis of the linearized KW operator}\label{analKW}
We now present analytic methods for the study of the linear and nonlinear KW operators $\LKW$ and ${\KW}$.  The goals of this analysis
have already been used above: first to prove that $\LKW$ is Fredholm between certain natural function spaces, and to calculate its index; 
and second, to show that solutions of ${\KW}(A, \phi)$ are polyhomogeneous near the knot $K$, justifying the calculations in the Weitzenb\"ock 
formula and the uniqueness theorem. These are direct generalizations of the results concerning the case $K = \emptyset$ in 
\cite{MW}, but the proofs here are complicated by the more singular nature of the problem. 

Let $(M,g)$ be a four dimensional compact Riemannian manifold with boundary $W$, $E$ an $\mathrm{SU(2)}$ or $\mathrm{SO(3)}$ bundle 
over $M$, and $K \subset W$ a closed knot or link. Choose local coordinates $\vec x = (x^1, x^2, x^3, x^4)$ as above, where $x^1 = t$ is 
arclength along $K$, $x^2, x^3$ are Fermi coordinates for the restriction of $g$ to $W$ around $K$ in $W$, and $y = x^4$ is geodesic distance 
from $W$, all with respect to $g$.   The corresponding cylindrical coordinates around $K$ are $\rho = |(y, x^2, x^3)|$ and 
$\omega = (\omega_0, \omega') = (y, x^2, x^3)/\rho \in S^2_+ = (\sin s, \cos s \cos \theta, \cos s \sin \theta)$, $s \in [0,\pi/2]$, $\theta \in [0,2\pi]$. 
The background metric
\begin{equation}
g = \d \rho^2 + \rho^2( \d s^2 + \cos^2 s \d \theta^2) + \d t^2 + \ \mbox{higher order terms.} 
\label{g2}
\end{equation}
is incomplete. Certain constructions below are phrased in terms of the conformally related complete metric
\begin{equation}
\hat{g} = \frac{\d \rho^2 + \d t^2}{\rho^2 s^2} +  \frac{\d s^2 + \cos^2 s \d \theta^2}{s^2} + \ \mbox{higher order terms.} 
\label{ghat}
\end{equation}
The higher order terms in each case refer to tensors which decay relative to the displayed leading part of each of
these metrics. 

Let $(A,\phi)$ be a solution of the gauged KW equations on $(M,K)$, with a Nahm pole singularity along $W$ away from $K$ and a knot singularity
along $K$.  Denote by $\LKW$ the linearization of $\KW$ at this solution, acting on the infinitesimal variations $(a,\varphi)$. 
We repeat for convenience that 
\begin{equation}
\begin{split}
\LKW: \calC^\infty( & M, (\Lambda^1 \oplus \Lambda^3) \otimes \ad \, \mathfrak g) \longrightarrow
\calC^\infty(M, (\Lambda^0 \oplus \Lambda^2 \oplus \Lambda^4) \otimes \ad \, \mathfrak g), \\[0.5ex]
& \LKW(a, \star \varphi) = \big(d_A + d_A^* + M_{\phi}\big) (a, \star \varphi),
\end{split}
\label{linLKW22}
\end{equation}
where
$$
M_{\phi}( a, \star \varphi) = \big( \star [ \phi, \star \varphi],  [\phi, \star (\star \varphi)] + \star [a, \phi], [a, \star \phi]\big). 
$$
Near $W$ but away from $K$, this has the form 
$$
\LKW = \d + \d^* + \frac{1}{y} B_0,
$$
or more pertinentily here, near $K$ 
$$
\LKW = \d + \d^* + \frac{1}{\rho s} B_0,
$$
where $B_0$ is an endomorphism bounded up to $\rho = 0$ and $s=0$ (so also $y=0$ away from $K$). 
It is also useful to consider the operator
\begin{equation}
\hat{\LKW} = \rho s (\d + \d^*) + B_0,
\label{2LKW}
\end{equation}
which stands in relation to $\LKW$ as the complete metric $\hat{g}$ does to the incomplete metric $g$. (More specifically,
it is a complete iterated edge operator; however, the operator $\d^*$ here is still with respect to $g$, so this is not the
linearized KW operator for $\hat{g}$.)

We focus on the analytic properties of $\LKW$ near $K$, since its behavior near $W\setminus K$ has already
been treated in \cite{MW}.  This reduction depends on the fact that the results there are local in $W$ so may be transported
to the present setting. 

As in \cite{MW}, the main part of this analysis is the construction of a parametrix $G$ for $\LKW$ near $K$. We then investigate
various properties of $G$ and use these to deduce our main results.  These steps are close to those in our earlier paper,
but with the additional difficulties caused by the more intricate singular structure of the operator. 

The reader may find the following guide to the rest of the paper useful.   In the remainder of this section we describe how this analysis
fits into a larger and more systematic framework from geometric microlocal analysis.  This entails the introduction of a manifold with corners 
$M_K$, the blowup of $M$ around $K$, and a brief explanation of a general inductive strategy to analyze the class of elliptic ``incomplete
iterated edge operators''.  The operator encountered in \cite{MW} occupies the first step in this inductive scheme, while the operator
here is of ``depth two''.   There is a whole category of objects associated to this class of operators, including the manifold with 
corners $M_K$ and various generalizations of it, Lie algebras of vector fields, classes of symbols, and hierarchies of model operators.
We give here a fairly minimal treatment, presenting only what is needed for the immediate purposes.  This general
analytic program has its foundations, in some sense, in \cite{Ma}, with some parts of the inductive strategy laid out in \cite{ALMP, ALMP2}.  
The results needed for the applications here go beyond what is proved in these last two papers. 

After explaining this general program, we take up the parametrix construction. For expository purposes we first describe a rather
crude $L^2$ parametrix for $\LKW$. This brings out some of the essential new features of the problem, but we cannot deduce
enough from it to deduce, for example, the mapping properties of $\LKW$ on weighted H\"older spaces.  We therefore present
the more involved construction of the geometric microlocal parametrix.  From this we are able to prove various refined mapping
properties of $\LKW$.  

The final section of this paper contains the proofs of the regularity theorems. These are based on commutator arguments with the parametrix
and a new iterative scheme to improve regularity.  

\subsection{The space $M_K$}
The resolution $M_K$ is the radial blowup of $M$ around $K$ (in general such a blowup is denoted $[M;K]$), i.e., the disjoint union
of $M \setminus K$ and the inward-pointing spherical normal bundle of $K$. Its structure as a manifold with corners is reflected 
by the cylindrical coordinates $(\rho, t, \omega)$ around $K$, where $\omega \in S^2_+$, $\rho \geq 0$ and $t$ is a local 
coordinate along $K$, and indeed these are a nondegenerate coordinate system on $M_K$. We also use coordinates $\omega = (s,\theta)$
on the hemispheres $S^2_+$, where $s=0$ is the boundary and $\theta \in S^1$. Thus $M_K$ is endowed with the unique minimal 
smooth structure generated by the lifts of smooth functions on $M$ and these cylindrical coordinates. 

This blowup is a manifold with corners of codimension $2$ and has two boundary hypersurfaces. The first, denoted $\wfa$,  is the 
closure of the lift of $W \setminus K$; the second is the `new' boundary hypersurface $\ff$ created by the blowup and called 
the front face of $M_K$. The function $\rho$ is a boundary defining function for $\ff$, so $\ff = \{\rho = 0\}$ and $\d \rho \neq 0$ 
there, while $s$ is a boundary defining function for $\wfa$.   Notice that $\ff$ is the total space of a fibration over $K$, where 
each fibre is a copy of $S^2_+$, which we sometimes also write as $Z$. 

The passage to $M_K$ focuses attention on the approximate homogeneities of the analytic problem. For example, under favorable conditions, 
solutions of $\LKW u = 0$ are polyhomogeneous on $M_K$; viewed just on $M$, however, their singularities at $K$ look much
less tractable. Similarly, it is straightforward to check that the lift to $M_K$ of the approximate solution $(A^K,\phi^K)$ is also polyhomogeneous. 

\subsection{Stratified spaces and iterated edge operators}
We now briefly review some points about stratified spaces, referring to \cite{ALMP} for a careful account of the structural 
axioms and geometry. The key is the inductive way in which these spaces are constructed and many facts about them 
proved. The simplest class of stratified spaces are the smooth compact manifolds; these have only one stratum,
and by definition have depth $0$. There is an operation which constructs a new stratified space out of an old one,
increasing the depth (or complexity). Namely, if $Z$ is any compact stratified space of depth $k$, we may pass to 
a truncated cone over $Z$, or slightly more generally, a bundle of such cones.  These are the models for neighborhoods 
of the most singular strata in a space with depth $k+1$.   Thus a compact stratified space of depth $1$ decomposes as 
the union of a smooth open manifold and an open truncated cone or bundle of cones over a compact, smooth (depth $0$) space. 
Depth $1$ spaces are those with isolated conic singularities or, if the singular set has dimension greater than zero, simple 
edge singularities.  Continuing in the same way, a compact space of depth $k+1$ decomposes into an an open space which
has only depth $k$ singularities and another open space which is a tubular neighorhood of the depth $k+1$ stratum, and which 
is a cone or bundles of cones over a compact smooth base and with cross-section a compact depth $k$ space.

There is an important resolution procedure: successively blowing up the singular strata of a depth $k$ stratified space $X$ 
(in order of decreasing depth of the strata) yields a compact manifold with corners up to codimension $k$, $\wh X$.
The interiors of $X$ and $\wh{X}$ are canonically identified. The boundary hypersurfaces of $\wh X$ correspond 
to the blowups of each of the singular strata of $X$. Each of these boundary hypersurfaces is the total spaces of a 
fibration, where the fibers are the resolutions of the cross-sections of the cone-bundles associated to that stratum. These 
fibrations fit together in a precise way at the corners of $\hat{X}$. 

Let us specialize immediately to our case of interest. Decompose $M$ into strata $(M \setminus W) \sqcup 
(W \setminus K) \sqcup K$. As described above, a neighborhood  $\calU$ of $K$ in $M$ is diffeomorphic to a 
bundle of cones over $K$, with cross-section the hemisphere $S^2_+$. 
We often write $S^2_+= Z$ below. The resolution $M_K$ is obtained by replacing each point $p \in K$ by its inward-pointing 
spherical normal bundle, which is a copy of $S^2_+$.  
The two boundary hypersurfaces of $M_K$, $\wfa$ and $\ff$, correspond to the original boundary $W \setminus K$ and 
the front face produced in the blowup, respectively. 

On any stratified space there is a distinguished class of incomplete metrics and the associated differential operators,
called $\iie$ metrics and $\iie$ operators. We describe these only on $M$. Using the coordinates $(\rho, t, s, \theta)$,
an $\iie$ operator $B$ of order $1$ takes the form
$$
B = B_\rho \del_\rho + B_t \del_t + \frac{1}{\rho}\left( B_s \del_s + B_\theta \del_\theta  + \frac{1}{s}B_0 \right), 
$$
where each $B_j$ is an endomorphism smooth (or polyhomogeneous) on all variables. The operator in parentheses, 
$$
B_Z := B_s \del_s +  B_\theta \del_\theta + \frac{1}{s}B_0,
$$
is an $\iie$ operator on the depth $1$ space $Z$; in fact, it is just the two-dimensional version of 
the operator studied in our earlier paper \cite{MW}. 

Referring back to \eqref{linLKW22}, since $\phi \sim \rho^{-1} \phi_0$ and $\phi_0 \sim s^{-1}$, we see that the
linearized KW operator $\LKW$ is an $\iie$ operator.  For much of this section, we use little about the precise 
structure of its coefficients, only their asymptotic form near $\rho = 0$ and $s=0$, though certain things
(e.g., the computation of indicial roots) require more information. Thus we write
\begin{equation}
\LKW = B_\rho \del_\rho + B_t \del_t + \frac{1}{\rho} \mathcal J,
\label{linKW}
\end{equation}
where 
\begin{equation}
\calJ = B_s \del_s +  B_\theta \del_\theta  + \frac{1}{s}B_0
\label{linKW2}
\end{equation}
is a uniformly degenerate operator on $Z$.  

An $\iie$ operator $B$ is called $\iie$-elliptic if its $\iie$ symbol is invertible. In more detail, set
$$
\rho s B = B_\rho (\rho s \del_\rho) + B_t (\rho s \del_t) +  B_s (s \del_s) + B_\theta s\del_\theta + B_0,
$$
formally replace $\rho s \del_\rho$ by $-i\xi$, $\rho s \del_t$ by $-i\tau$, $s\del_s$ by $-i\sigma$ and $s\del_\theta$ by $-i\zeta$,
and drop the lower order term $B_0$. This defines the $\iie$ symbol
\begin{equation}
{}^\iie\sigma_1(B)(\rho, t, s, \theta; \xi, \tau, \sigma, \zeta) = \frac{1}{i}\left(B_\rho  \xi + B_t  \tau + B_s \sigma + B_\theta \zeta \right).
\label{iiesymbol}
\end{equation}
(There is a way to interpret this invariantly, but we do not do this here.) We require that this endomorphism
is invertible whenever $(\xi, \tau, \sigma, \zeta) \neq 0$. It is straightforward to calculate that $\LKW$ is $\iie$ elliptic.

\section{Parametrices}\label{param}
The basic problems about an elliptic $\iie$ operator are to show that it is Fredholm acting between appropriate 
function spaces, to compute its index, and to analyze the regularity of solutions of $B u = f$. These will be 
addressed using parametrix methods.  

\subsection{The parametrix construction in the simple edge case} 
To set the stage, we review how this works in the simple edge case, cf.\ the lengthier discussion in \cite{MW}. To be concrete, 
let us consider a problem very similar to ours and assume that $M$ is a closed $4$-dimensional manifold and $K$ 
an embedded knot. A tubular neighborhood of $K$ in $M$ is diffeomorphic to a bundle of truncated  cones $C_1(Z)$, where now 
$Z = S^2$. The blowup of $M$ around $K$ is the space $M_K = [M;K]$, which is a manifold with boundary $\del M_K$. This boundary
is the total space of a fibration over $K$ with fiber $Z$. Cylindrical coordinates around $K$ lift to a nondegenerate local coordinate system 
$(\rho, t, z)$, where $\rho$ is a boundary defining function for $\del M_K$ which restricts to a radial function on each conic fiber, 
while $t$ and $z$ are local coordinates on $K$ and $Z$, respectively.  The lift of a Riemannian metric $g$ on $M$ which is smooth
across $K$ takes the form
\begin{equation}
g = g_0 + \eta, \quad \mbox{where}\  g_0 = \d \rho^2 + h_K + \rho^2 k_Z;
\label{edgemetric}
\end{equation}
here $h_K$ is the restriction of $g$ to $K$ at $\rho=0$, $k_Z$ is the spherical metric on $Z$ (both are extended to this tubular 
neighborhood), and the tensor $\eta$ is a lower order deviation in the sense that $|\eta|_{g_0} \leq C \rho^\epsilon$ 
for some $\epsilon > 0$ as $\rho \searrow 0$. More generally, of course, one could consider spaces where $Z$ is not a sphere so
the locus $K$ is a genuine topological singularity.   Note the similarity with the metric \eqref{g2} -- indeed, the only difference
is that here $Z$ is a closed manifold while in \eqref{g2} it is a manifold with boundary. 

Consider a general Dirac-type operator $B$ adapted to the pair $(M,K)$, so for example $B$ is the sum
of an operator which is smooth across $K$ and a term $B_0/\rho$ of order $0$ which blows up like $\rho^{-1}$. We suppose
that all coefficients are smooth in cylindrical coordinates. We are tacitly assuming that $B$ satisfies an edge symbol ellipticity 
condition, see below for the corresponding condition in the depth $2$ setting. The indicial data of $B$ is the set of pairs $(\gamma, \Psi_0)$ 
where $\Psi_0$ is a field on $Z$ such that $B( \rho^\gamma \Psi) = \calO(\rho^{\gamma})$ (rather than the expected order 
$\calO(\rho^{\gamma-1})$ for some smooth extension $\Psi$ of $\Psi_0$.  As before, this means that $\gamma$ solves 
a generalized eigenvalue problem on $Z$ and $\Psi_0$ is the corresponding eigenfunction.  

Indicial roots determine some features of the mapping properties of $B$. The other ingredient is the normal operator $N(B)$:
for any $t_0 \in K$, $N(B)$ is the scale- and translation-invariant operator on $\RR \times C(Z)$ (which in our case is just 
$\RR^4$ since $Z = S^2$) modeling $B$ at that point. It is obtained by freezing the coefficients of $B$ at $\rho=0$, $t=t_0$
and introducing a certain global set of projective coordinates $(r, \tau, \omega)$ on $\RR^+ \times \RR \times S^2$. 
The full rationale for these new global coordinates will be explained later, but for the moment we regard this as a formal
change. In terms of these, we set 
\begin{equation}
N(B) = B_\rho(0,t_0,z) \del_r + B_t(0,t_0,z) \del_\tau + \frac{1}{r} \calJ,
\label{rhonorm}
\end{equation}
where $\calJ$ is a Dirac-type operator on $Z$ and $B_\rho(0,t_0)$, $B_t(0,t_0)$ are endomorphisms.  These new coordinates
reduce two important features: $N(B)$ is translation invariant in $\tau$ and homogeneous with respect to dilations 
$(r,\tau) \mapsto (\lambda r, \lambda \tau)$.   For more general $\iie$ operators, it is possible for both the indicial data and
the normal operator $N(B)$ to depend nontrivially on $t_0 \in K$; for simplicity we assume that this is not the case, and
fortunately our particular problem does not have this complicating feature. 

We next define weighted edge Sobolev spaces on $M_K$: 
\begin{definition} The space $\rho^{\mu + 3/2} H^\ell_e (M_K)$ consists of the space of functions (or fields) $u = \rho^{\mu + 3/2} v$
where $(\rho\del_\rho)^a (\rho\del_t)^b \del_z^c v \in L^2$ for all $a + b+ c \leq \ell$. 
\end{definition}
The vector fields appearing here generate the bounded vector fields for the metric $\hat{g} = \rho^{-2} g$, and hence $H^\ell_e(M_K)$ is
the natural Sobolev space associated to that complete metric. The shift by $3/2$ in the weight parameter is a normalization:
if $\Psi\in \rho^{\mu+3/2} L^2$ is supported in $\{\rho \leq 1\}$, then since $dV_g \cong \rho^2 \d \rho d\tau dz$, 
$\int \rho^{-2(\mu+ 3/2)} |\Psi|^2 \, \rho^2 \d \rho d\tau dz = \int |\Psi|^2 \rho^{-2\mu-1}\, \d \rho d\tau dz< \infty$ if $|\Psi| \leq C 
\rho^{\mu + \epsilon}$ for any $\epsilon > 0$.  In other words, with this normalization, then $\Psi \precsim \, \rho^\mu$ in an $L^2$ sense. 

We also define $r^{\mu + 3/2} H^\ell_e(C(Z) \times \RR)$, where now all integrals are taken on the entire product space. 
We say that $\mu$ is an elliptic weight for $B$ if
$$
N(B): r^{\mu + 3/2} H^1_e( C(Z) \times \RR) \longrightarrow r^{\mu  + 1/2} L^2 (C(Z) \times \RR)
$$
is an isomorphism. 

There are corresponding edge Sobolev spaces spaces on $M_K$ defined with respect to differentiations by the vector 
fields $\rho\del_\rho$, $\rho \del_t$, $\del_z$.  \cite[Theorem 6.1]{Ma} states that if $\mu$ is an elliptic weight for $B$, then 
\begin{equation}
B:  \rho^{\mu  + 3/2} H^{k+1}_e(M_K) \longrightarrow \rho^{\mu+1/2} H^k_e(M_K)
\label{fredBSob}
\end{equation}
is Fredholm for every $k$.  We may also define weighted H\"older spaces $\rho^\mu \calC^{k,\alpha}_e$ relative to 
differentiations by the same set of degenerate vector fields, and by \cite[Corollary 6.4]{Ma}, 
\begin{equation}
B: \rho^\mu \calC^{k,\alpha}_e(M_K) \longrightarrow \rho^{\mu-1} \calC^{k-1, \alpha}_e(M_K)
\label{fredBHold}
\end{equation}
is also Fredholm.  We emphasize that $B$ is not Fredholm for more general weights, either because it does not have closed range or
else because it has infinite dimensional kernel or cokernel.

These results are proved by first constructing an approximate inverse, or parametrix, to $B$, i.e., an operator 
$G = G_\mu$ which satisfies 
$$
B \circ G = \mbox{Id} + Q_1, \qquad  G \circ B = \mbox{Id} + Q_2,
$$
where $Q_1$ and $Q_2$ are compact on $\rho^{\mu}H^k_e$ and $\rho^{\mu+1}H^{k+1}_e$, respectively, and
in addition is such that 
$$
G: \rho^{\mu} H^k_e(M_K) \longrightarrow \rho^{\mu + 1} H^{k+1}_e(M_K)
$$
is bounded. This implies that \eqref{fredBSob} is Fredholm.  Showing the analogous property for \eqref{fredBHold} relies on a 
description of the pointwise behavior of the Schwartz kernel of $G$, from which it can be proved that
$$
G: \rho^{\mu-1}\calC^{k, \alpha}_e(M_K) \longrightarrow \rho^{\mu} \calC^{k+1,\alpha}_e(M_K)
$$
is bounded, which in turn shows that \eqref{fredBHold} is Fredholm. 

This `sharp structure' of $G$ is described as saying that its Schwartz kernel is polyhomogeneous on a certain further resolution 
of $M_K^2$, apart from a `uniform' singularity along the lifted diagonal.  We explain this. Recall that $G$ is a distribution on $M_K^2$ 
which is smooth in the interior away from the diagonal. It has a standard pseudodifferential singularity across the diagonal, but is 
also singular along $\del M_K \times M_K$ and $M_K \times \del M_K$. Its most `important' singularity occurs at the fiber diagonal 
of $(\del M_K)^2$, and to understand it we resolve $M_K^2$ by blowing up along the fiber diagonal to obtain a new and slightly more 
complicated manifold with corners $(M_K)^2_e$. This is a manifold with corners of codimension $2$, with new boundary 
hypersurfaces called the front face. The lift of the diagonal of $M_K^2$ to this space intersects this front face transversely.

The paper \cite{Ma} considers the class of edge pseudodifferential operators, which are by definition those operators
whose Schwartz kernels are polyhomogeneous on $(M_K)^2_e$, along with a classical singularity along the lifted diagonal. 
These are modeled on the type of degenerate behavior already seen in the differential operator $B$.
The technical results there involve understanding the compositions of such operators and the boundedness of these
operators on weighted Sobolev and H\"older spaces. A central conclusion is that if $\mu$ is an elliptic weight for $B$,
then one can construct a parametrix $G$ for $B$ in this space of pseudodifferential operators with compact
remainders $Q_j$. The general boundedness properties of these operators lead to the assertions above.  

To expand on this slightly, let $\mathrm{fdiag}$ be the fiber diagonal in $(\del M_K)^2$. Using $(\rho, t, z, \tilde{\rho},
\tilde{t}, \tilde{z})$ as local coordinates on $(M_K)^2$, this fiber diagonal is $\{\rho = \tilde{\rho} = 0, t = \tilde{t} \}$.
The edge double space $(M_K)^2_e$ is the blowup $[ M_K^2; \mathrm{fdiag}]$.  Edge pseudodifferential operators, and
in particular, the parametrix $G$, are polyhomogeneous at the various boundaries of this space. Thus, beyond
its singularity across the lifted diagonal, $G$ has complete asymptotic expansions at each of the boundary hypersurfaces 
of $(M_K)^2_e$ and product-type expansions at the corners. The exponents in these expansions at the left and right `side' 
faces corresponding to $\del M_K \times M_K$ and $M_K \times \del M_K$ are determined by the indicial roots of $B$; 
the leading exponent at the front face is the universal number $1 - \dim M = -3$.  The boundedness properties
of $G$ are determined  by the leading terms of these expansions.  The parametrix also leads to higher regularity
results, including results like if $B \Psi = 0$ and $\Psi \in \rho^{\mu +1}H^1_e$ where $\mu$ is an elliptic weight,
then $\Psi$ is polyhomogeneous on $M_K$.

\subsection{The parametrix for $\LKW$ in depth $2$}
In our actual problem, $K$ lies on the boundary of $M$ and each conic cross-section $Z$ has a boundary, so 
the structure of a parametrix $G$ for $\LKW$ is more complicated near the boundaries and corners of $M^2_K$.  
We take a slightly circuitous route to describe this.  We first show that the same objects which appeared in the
simple edge theory above, indicial roots and the normal operator, can be defined here as well.  A more elementary
parametrix construction yields a rough parametrix, which implies that $\LKW$
is Fredholm on weighted Sobolev spaces.  We then carry out the geometric microlocal parametrix construction,
where the aforementioned `new' singularities of $G$ appear explicitly. This leads to the fine regularity properties 
for solutions in weighted H\"older spaces. 

\subsubsection{Indicial roots}
As in the simple edge case, the indicial roots of $\LKW$ are the formal rates of decay or growth of solutions of the linear problem 
$\LKW \Psi = 0$, and also play a role in formulating the global mapping properties of this operator. As we have described in
Sec.4, there are different indicial root sets at each of the two boundaries $\del_W M_K$ and $\del_K M_K$; those in the first case were already 
calculated in \cite{MW}, while the calculations of the ones at $\rho=0$ occupy the first part of this paper. 
In this and later sections we explain how they enter the analytic theory. In a certain sense the story is not so different from
the simple edge case once suitable boundary conditions for $\calJ$ are imposed at $\del Z$.

We say that $(\gamma, \Psi_0(z))$ is an indicial pair for $\LKW$ at $\del_K M_K$ (or a $\rho$-indicial pair) if there is a smooth field 
$\Psi \sim \Psi_0 + \rho \Psi_1 + \ldots$ such that 
$$
\LKW (\rho^\gamma \Psi)= \mathcal O(\rho^\gamma). 
$$
The precise extension $\Psi$ of $\Psi_0$ is unimportant. By \eqref{linKW}, 
\begin{equation}
\LKW (\rho^\gamma \Psi) = \rho^{\gamma-1} \left( B_\rho(0,t,z) \gamma + \left.\calJ \right|_{\rho=0} \right) \Psi_0 + 
\mathcal O(\rho^\gamma),
\label{indface}
\end{equation}
so $\gamma$ is an indicial root if $\left. \left(  \calJ + B_\rho \gamma \right)\right|_{\rho=0}$ has nontrivial nullspace, 
and then $\Psi_0$ must lie in this nullspace. The variable $t \in K$ enters as a parameter; in principle the indicial roots 
might depend on $t$, but fortunately that is not the case here.

As noted above, there are also indicial roots for $\LKW$ at $\del_W M_K$, which we call the $s$-indicial roots for short.  These were
computed in \cite{MW}.  It is worth noting that unlike the $\rho$-indicial roots, these are determined by an algebraic 
eigenvalue equation. 

When $Z$ is closed, the set of $\rho$-indicial roots is always discrete, but when $Z$ is a manifold with boundary, we must 
impose boundary conditions for $\calJ$ at $\del Z$ to obtain a self-adjoint problem. (In the absence of such boundary 
conditions, solutions $\Psi_0$ would exist for any value of $\gamma$.)  In the language of \cite{MW}, we are in the 
quasi-regular case, which means that $0$ is not an $s$-indicial root. It therefore suffices let $\cal$ act on fields which
vanish at $s=0$, i.e., which satisfy Dirichlet boundary conditions at $\del Z$. 

\subsubsection{$L^2$ Sobolev spaces}
Consider first the weighted $L^2$ spaces $\rho^\mu s^\nu L^2(M_K, dV_g)$ for any $\mu, \nu \in \RR$, where up to
a bounded nonvanishing factor, $dV_g = \rho^2 \d \rho \d t \d s \d \theta$. We use Sobolev spaces adapted 
to the degeneracy structure of the problem encoded by the vector fields $\rho s \del_\rho$, $\rho s \del_t$, $s\del_s$, $s\del_\theta$.
Omitting bundles from the notation notation unless necessary, define 
$$
\rho^\mu s^\nu H^k_{\ie}(M_K, dV_g) = \{ \Psi:  (\rho s \del_\rho)^i (\rho s \del_t)^j (s\del_s)^\alpha (s\del_\theta)^\beta \Psi \in 
\rho^\mu s^\nu L^2(M_K), \ \ \forall \ \ i + j + \alpha + \beta \leq k\}
$$
Clearly, since $\rho s \LKW$ is a smooth combination of a generating set for $\calV_\ie$, 
$$
\LKW:  \rho^{\mu} s^{\nu } H^{k+1}_{\ie}(M_K) \longrightarrow \rho^{\mu  -1 }s^{\nu-1} H^k_{\ie}(M_K)
$$
is bounded for every $\mu, \nu$.  

Recalling the quasi-regularity assumption, let $(\underline{\lambda}, \overline{\lambda})$ be the maximal interval containing $0$ 
which contains no $s$-indicial roots.  We explain later that nontrivial solutions to the generalized eigenvalue problem 
$(\calJ + \gamma B_\rho)\Psi_0 = 0$ with $\Psi_0 \in s^{\lambda+1/2} L^2(Z, \d s\d \theta)$ for any $\lambda$ in this interval exist 
only when $\gamma$ is one of the $\rho$-indicial roots computed in Section \ref{indicial}.

Define $\calD$ as the space of all fields $\Psi_0$ on $Z$ such that $\Psi_0 \in s^{\lambda+1/2}L^2$ and $\calJ \Psi_0 \in s^{\lambda + 1/2}$.  Now set
\begin{equation}
\rho^{\mu + 3/2} H^k_{\ie, \calD}(M_K) = \{\Psi \in \rho^{\mu + 3/2} H^k_{\ie}:  \Psi(\rho, t, \cdot) \in \calD\ 
\mbox{for a.e.} \ \rho, t \}
\label{defieD}
\end{equation}
with norm 
$$
||\Psi||_{\ie; k, \mu + 3/2, \calD}^2 = || \Psi||^2_{\ie, k, \mu + 3/2, \lambda + 1/2} + 
|| \calJ \Psi||^2_{\ie, k-1, \mu + 3/2, \lambda + 1/2}.
$$

We shall eventually prove that 
\begin{equation}
\LKW: \rho^{\mu + 3/2} H^k_{\ie, \calD}(M_K) \longrightarrow \rho^{\mu + 1/2} H^{k-1}_{\ie}(M_K).
\label{mapwl2}
\end{equation}
is Fredholm if and only if $\mu$ is not a $\rho$-indicial root of $\LKW$. 

\subsubsection{Ellipticity and the normal operator}
As in the simple edge case, three hypotheses constitute the definition of full ellipticity. The first is the invertibility of the $\ie$ symbol
discussed earlier, while the second and third involve the invertibility of the $s$- and $\rho$-normal operators.  These normal
operators are the models for $\LKW$ at $\del_W M_K$ and $\del_K M_K$, respectively. We carry over
from \cite{MW} the fact that in the range of $s$-weights $\lambda \in (\underline{\lambda}, \overline{\lambda})$, the
$s$-normal operator is invertible.  The remaining condition involves the mapping properties of the $\rho$-normal operator. 
Recall from \eqref{rhonorm} that this is the operator, parametrized by $t_0 \in K$, and acting on the entire space $\RR^2_+ 
\times S^2_+$, given by 
\begin{equation}
N_\rho(\LKW) =  B_\rho(0,t_0) \del_r + B_t(0,t_0) \del_t + \frac{1}{r} \calJ.
\label{normalop}
\end{equation}
Here $B_\rho$, $B_t$ and the coefficients of $\calJ$ are evaluated at $\rho=0$, $t = t_0$.
Clearly $N_\rho(\LKW)$ is the linearized $\KW$ operator on $(\RR^4_+)_{\RR} = [ \RR^4_+;  \RR]$ in cylindrical coordinates,
and is independent of $t_0$. 

The value $\mu$ is called an elliptic weight (in $\rho$) if 
\begin{equation}
N_\rho(\LKW):  \rho^{\mu + 3/2} H^1_{\ie, \calD}((\RR^4_+)_{\RR}) \longrightarrow \rho^{\mu+ 1/2} L^2((\RR^4_+)_{\RR})
\label{normalellip}
\end{equation}
is invertible.  The shift in weight parameter is a normalization which reflects that $\rho^\gamma \in \rho^{\mu + 3/2} L^2 (\rho^2 \d \rho)$ 
(locally near $\rho=0$) if and only if $\gamma > \mu$. 

In view of the formulation of the Nahm pole boundary condition around a knot, we choose $\mu \in (-1,0)$ to
allow perturbations of the approximate solution which grow at most like $\rho^{\mu}$, i.e., like $\rho^{-1 + \varepsilon}$. 
The indicial root computations of section 4, see in particular Proposition \ref{fullindrts}, shows that this interval is 
free of indicial roots, so \eqref{normalellip} has closed range.  
\begin{prop}
If $\mu \in (-1,0)$, the mapping \eqref{normalellip} is an isomorphism. 
\label{normalisom}
\end{prop}
The proof relies on two facts: this map is injective, and has index zero.  

The Fredholm and regularity theory proved below show that neither of these properties change so long as $\mu$ 
does not cross an indicial root, so we may as well set $\mu = -1/2$. We then have that $\rho H^1_{\ie, \calD}$
is precisely the set of $u \in L^2$ such that $\LKW u \in L^2$, or in other words, $\LKW: \rho H^1_{\ie, \calD} \to L^2$ 
is simply the realization of the unbounded operator $\LKW: L^2 \to L^2$ to its domain.  The discussion in section 6
shows that $\LKW: L^2 \to L^2$ is pseudo-skew adjoint; the noncompactness of the `cross-section' does not affect this,
so we conclude that its index vanishes, as required.

Now suppose that $N_\rho(\LKW) \Psi = 0$, where $\Psi \in \rho L^2$.  Recall that $N_\rho(\LKW)$ is just
the linearization $\LKW$ of the KW operator at the model knot solution, and following the discussion
in section 5, $\LKW^\dag \LKW = \hD$, where $\hD$ is the operator in \eqref{NSPS}. We showed there
that we may decompose into eigenmodes, first for the induced operator on $S^2_+$, so the coefficients 
satisfy \eqref{eigeneqn}, and then further into eigenmodes for the equation in $\omega \in S^1_+$. 
The coefficients now satisfy a homogeneous Euler equation, and hence are monomials $R^\sigma$.
The injectivity now follows, since any such monomial which is in $L^2$ near $R=0$ fails
to be in $L^2$ as $R \to \infty$, and vice versa.   As explained above, we deduce from this that
$\LKW$ is also surjective on $L^2$.  This is one of the key facts needed in the geometric microlocal parametrix
construction. 

\subsubsection{The $L^2$ parametrix}
Writing $\LKW$ as in \eqref{linKW}, the coefficient endomorphisms $B_\rho$, $B_t$ and the operator $\calJ$ all depend smoothly 
on $t$, and hence vary slowly in this variable. Cover $K$ by neighbourhoods in $M$ which have the product form 
$\calU_\ell \cong [0, \rho_0)_\rho \times (t_\ell -\epsilon, t_\ell + \epsilon)_t \times Z$. For each such neighborhood, denote by $\LKW_\ell$ 
the operator $\LKW$ with coefficients frozen at $\rho = 0$, $t = t_\ell$.  For each $\ell$, write $\calE_\ell = \LKW - \LKW_\ell$..
Adjoin also the open set $\calU_0 = M \setminus \{ \rho \leq \rho_0/2\}$. The error  term vanishes at $\{(0, t_\ell\} \times Z$,
and hence is small throughout $\calU_\ell$.  Furthermore, $\LKW_\ell$ equals the normal operator $N(\LKW)$ at $(0,t_\ell)$.
In the first approach, we analyze it using the eigendecomposition for $\calJ$. 

An important subtlety is that the domain $\calD$ for $\calJ$ does not depend on the basepoint $t \in K$. This is because
it is defined by weight conditions involving the indicial roots for $\calJ$. 

Denote by $\LKW_{\ell,\lambda}$ the reduction of $\LKW_\ell$ to the eigenspace associated to any eigenmode $\lambda$ of $\calJ$,  
and write the associated eigencomponent of $\Psi$ as $\Psi_\ell$. Then 
$$
\rho \LKW_{\ell,\lambda} = B_\rho \, \rho \del_\rho + B_t \rho \del_t + \lambda.
$$
is a uniformly degenerate operator. 

For simplicity, assume that $t_\ell = 0$ and drop the subscript $\ell$ from the notation. Set  $\Sigma = [0,\rho_0) \times (-\epsilon, 
\epsilon)$. Uniformly degenerate differential operators on $\Sigma$ are sums of products of the basic vector fields $\rho \del_\rho$ 
and $\rho \del_t$. As described earlier and more fully in \cite{MW}, the corresponding class of pseudodifferential operators 
$\Psi_0^*(\Sigma)$ contains operators $H$ for which the Schwartz kernels lift to a polyhomogeneous distribution on the $0$ double 
space $\Sigma^2_0 = [\Sigma^2; \mathrm{diag}(\del \Sigma)^2]$. Coordinates on this space are given by the polar variables 
$R = |(\rho, \tilde{\rho}, t-\tilde{t})|$, $\Omega = (\rho, \tilde{\rho}, t-\tilde{t})/R$ along with $\tilde{t}$. 
The requirement is that $H(R, \Omega, \tilde{t})$ has a standard pseudodifferential singularity along the lifted
diagonal $\{ \Omega = (1/\sqrt{2}, 1/\sqrt{2}, 0)\}$ uniformly up to the front face of $\Sigma^2_0$ and has 
polyhomogeneous expansions at each boundary faces and product type expansions at the corners.

As quoted above, if $L$ is a uniformly degenerate edge operator of order $1$ which is fully elliptic at some weight $\mu$, i.e., its $0$-symbol is 
invertible and $\mu$ is an elliptic weight, so that the normal operator is an isomorphism $r^{\mu+1/2} H^1_0(\RR^2_+) 
\to r^{\mu-1/2} L^2(\RR^2_+)$, then there is a parametrix $G \in \Psi^{-1, E}_0(\Sigma)$ such that $LG = I + Q_1$, $GL = I + Q_2$ 
for some operators $Q_j \in \Psi^{-\infty, E}_0(\Sigma)$, and 
$$
G: \rho^{\mu + 1/2} L^2( \Sigma) \longrightarrow \rho^{\mu + 3/2} H^1_{0}(\Sigma),
$$
is bounded.  The $L^2$ spaces here are taken with respect to the measure $\rho^2 \d \rho \d t$. 
The index set $E$ in this notation describes the exponents in the expansions of these Schwartz 
at the various boundary faces of $\Sigma^2_0)$. We do not assert that this map is Fredholm because we are
considering it locally in $\Sigma$ (and have not specified boundary conditions at the other boundaries of
this region). 

We apply this to each reduced operator $\LKW_\lambda$, thus obtaining a sequence of operators $G_\lambda \in 
\Psi^{-1, E(\lambda)}_0(\Sigma)$ which are generalized inverses for $\rho \LKW_\lambda$.   We need to know slightly more 
before we can conclude that $\rho \LKW$ itself has a generalized inverse up to a compact error. Namely, we must 
show that $\LKW_\lambda$ is not 
an isomorphism for at most finitely many values of $\lambda$. This may be done directly, but we do not explain this 
since it follows from the other parametrix construction below.

We have now shown that for each $\ell$ and $\lambda$ there is a $0$-pseudodifferential operator $G_{\ell, \lambda}$ 
such that 
$$
G_{\ell,\lambda}:  \rho^{\mu + 1/2} L^2(\Sigma) \longrightarrow \rho^{\mu + 3/2} H^1_0 (\Sigma)
$$
is bounded and 
$$
\LKW_{\ell,\lambda} G_{\ell, \lambda} = \mathrm{Id} + Q_{\ell,\lambda},  \quad 
G_{\ell, \lambda} \LKW_{\ell,\lambda} = \mathrm{Id} + Q_{\ell,\lambda}',
$$
where $Q_{\ell,\lambda}, Q_{\ell, \lambda'}$ are residual, i.e., have Schwartz kernels which are smooth and vanish to infinite order 
as $\rho \to 0$ on $( [0,\rho_0) \times (-\epsilon, \epsilon))^2$. We also write
$$
G_\ell = \sum_\lambda  G_{\ell, \lambda} \Psi_\lambda(z) \Psi_\lambda(\tilde{z}).
$$

Choose a smooth partition of unity $\chi_\ell$ for the open cover $\{ \calU_0, \ldots, \calU_N\}$ of $M_K$ and another family of 
smooth functions $\tilde{\chi}_\ell$ with $\mathrm{supp} \, \tilde{\chi}_\ell \subset \calU_\ell$ and $\tilde{\chi}_\ell = 1$ on
$\mathrm{supp} \,\chi_\ell$.  Now define
\begin{equation}
\tilde{G} = \sum_\ell  \tilde{\chi}_\ell(\rho, t, z)  G_{\ell}( \rho, t, z, \tilde{\rho}, \tilde{t}, \tilde{z}) 
\chi_\ell( \tilde{\rho}, \tilde{t}, \tilde{z}). 
\label{firstparam}
\end{equation}
We calculate that
$$
\LKW \tilde{G} =  \mathrm{Id} + \sum_\ell  
\bigg( Q_\ell + [ \LKW, \tilde{\chi}_\ell] G_{\ell}  + \tilde{\chi}_\ell \calE_\ell G_\ell \bigg) \chi_\ell = \mathrm{Id} - \tilde{Q}, 
$$
with a similar expression for $\tilde{G} \LKW$.  

\begin{prop}
If $\mu$ is an elliptic weight for $\LKW$, then 
$$
\LKW:  \rho^{\mu  + 3/2} H^{k+1}_{\ie}(M_K) \longrightarrow \rho^{\mu + 1/2} H^{k}_{\ie}(M_K)
$$
is Fredholm.
\end{prop}
It suffices to prove that $\tilde{Q}$ and $\tilde{Q}'$ are sums of a compact operator and an operator with small norm 
on $\rho^{\mu + 1/2} H^{k}_{\ie}(M_K)$ and $\rho^{\mu + 3/2} H^{k+1}_{\iie,\calD}(M_K)$, respectively.  We check the three terms 
in this sum separately.   Clearly each $Q_\ell\chi_\ell$ is compact since its Schwartz kernel is smooth and vanishes to infinite 
order at the boundary. Next, the support of the multiplication operator $[\LKW, \tilde{\chi}_\ell]$ is disjoint from 
the support of $\chi_\ell$, so $[\LKW, \tilde{\chi}_\ell] G_\ell \chi_\ell$ has Schwartz kernel which is smooth in the interior 
of $(M_K)^2_0$, and vanishes near the front faces.  Finally, using that the coefficients of $\calE_\ell$ are small in $\calU_\ell$,
we see that $\tilde{\chi}_\ell\calE_\ell G_\ell$ is an operator of order $0$ multiplied by a very small coefficient, and hence
has small norm. 
\hfill $\Box$

\subsection{The geometric microlocal $\ie$ parametrix}
The more refined parametrix construction addresses the fact that the singularity structure of the generalized 
inverse to $\LKW$ is polyhomogeneous only when lifted to %
the $\ie$ double space $(M_K)^2_{\ie}$,  obtained from $M^2_K$ by an iterated blowup.  We now define this space and 
show that it provides the right setting to amalgamate the solution operators for various model problems to obtain
a parametrix for $\LKW$.  This procedure is much the same as in the simple edge setting, but now there are two front faces
and a correspondingly larger set of model problems.

\subsubsection{The structure vector fields} 
The degeneracy structure of operator $\LKW$ motivates the entire construction. The coordinate 
expression shows that $\rho s \LKW$ is a linear combination of the vector fields $\rho s \del_\rho$, $\rho s \del_t$, $s\del_s$ 
and $s\del_\theta$, with smooth (or polyhomogeneous) endomorphism valued coefficients. We define the structure vector fields for 
the problem to consist of the space of all smooth multiples of these basic vector fields, i.e., 
$$
\calV_{\ie} = \mbox{span}\,\{ \rho s \del_\rho, \rho s \del_t, s\del_s, s\del_\theta\}.
$$
This is a locally finitely generated Lie algebra and a $\calC^\infty(M_K)$ module, which is invariantly defined 
as the space of all smooth vector fields on $M_K$ which are unconstrained in the interior, vanish at $\del_W M_K$ and lie tangent 
to the fibers at $\del_K M_K$.  This is a so-called boundary fibration structure on $M_K$ which agrees with the uniformly degenerate 
structure at $\del_W M_K$ and with the edge structure at $\del_K M_K$, but these two structures interact at 
the corner $\del_W M_K \cap \del_K M_K$. 

\subsubsection{The $\ie$ double space}
The $\ie$ double space $(M_K)^2_\ie$ accommodates two objectives:  first, the lift of the diagonal $\diagiie$ is a 
`$p$-submanifold', i.e., it meets all faces and corners of $(M_K)^2_\ie$ cleanly (the diagonal in $M_K^2$ does
not have this property), and second, the lifts of the generating vector fields for $\calV_{\ie}$ are smooth on $(M_K)^2_\ie$ 
and span the normal bundle to $\diagiie$. We now define this resolution by iterated blowup and verify that these goals are 
satisfied.

The diagonal $\mathrm{diag}(M_K^2)$ has two boundary hypersurfaces, $\diag_s$ and $\diag_\rho$, which are its
intersections with the corners $(\del_W M_K)^2$ and $(\del_K M_K)^2$, respectively. The intersection $\diag_s \cap \diag_\rho$ 
is the diagonal of $(\del_W M_K \cap \del_K M_K)^2$.  Now consider the generating elements of $\calV_{\ie}$ lifted to the left factor 
of $M_K$ in $M_K^2$. These are nonvanishing away from the boundaries and span the normal bundle of $\diag$
in the interior. On the other hand, they vanish at $\diag_s$ and are tangent to the (left) hemisphere fibers in $(\del_K M_K)^2$ along 
$\diag_\rho$. 

We resolve these degeneracies by  blowing up two further important submanifolds: $\fdiag_\rho$, the fiber diagonal
of $(\del_K M_K)^2$, i.e., the set of points $(p, \tilde{p}) \in (\del_K M_K)^2$ such that $\pi(p) = \pi(\tilde{p})$
where $\pi: \del_K M_K \to K$ is the fibration, and $\diag_s$, the diagonal of $(\del_W M_K)^2$. The $\ie$ double space 
is obtained by blowing up these fiber diagonals in the order listed above: 
$$
(M_K)^2_{\ie} = [ M_K^2; \fdiag_\rho,  \diag_s].
$$
We denote by $\beta: (M_K)^2_\ie \to M_K^2$ the blowdown map. 
We describe the structure of this space in more detail. For simplicity, write $F_s = \del_W M_K$ and $F_\rho = \del_K M_K$,
so $M_K^2$ has four boundary hypersurfaces: 
$$
F_{\rho,\ell} = F_\rho \times M_K,\ F_{s,\ell} = F_s \times M_K,\ F_{\rho,r} = M_K \times F_\rho,\ F_{s,r} = M_K \times F_s.
$$
The subscripts $\ell$ and $r$ indicate that the face is on the left or right factor of $M_K$.  The two new front
faces created in the blowup are denoted $\ff_\rho$ and $\ff_s$.  We shall henceforth use the notation that
$R_\rho$ and $R_s$ are defining functions for these two faces; thus $\ff_\rho = \{R_\rho = 0\}$ and $d R_\rho \neq 0$
there, and similarly for $\ff_s$ and $R_s$. 

The face $\ff_s$ is the total space of a fibration $\pi_s$
over $\diag_s$;  since $\diag_s$ has codimension $5$, the fibers of $\ff_s$  are $4$-dimensional spherical orthants 
(with two coordinates nonnegative). On the other hand, the fibration structure of $\pi_\rho: \ff_\rho \to \fdiag_\rho$
degenerates over $\fdiag_\rho \cap \diag_s$. We see this as follows. After the first blowup, in the intermediate space 
$[ M_K^2; \fdiag_\rho]$, $\diag_s$, the front face created by the blowup is a bundle of two-dimensional spherical
orthants  (with two coordinates nonnegative) over the base $\fdiag_\rho$.  The submanifold $\diag_s$ intersects
each fiber over $\fdiag_\rho \cap \diag_s$ over a single point in the interior, so after the second blowup, the
preimage of any point in $\fdiag_\rho \cap \diag_s$ is the union of a $4$-dimensional spherical orthant and
a $2$-dimensional spherical orthant blown up around an interior point. On the other hand, the restriction of $\pi_\rho$ 
to the interior of $\ff_\rho$ is a true fibration over the interior of $\fdiag_\rho$. 

It is useful to see all of this in coordinates. Since we are working in a neighborhood of the diagonal of $M_K^2$, we
can use the same coordinate system on each factor of $M_K$, so $(\rho, t, s, \theta, \tilde{\rho}, \tilde{t}, \tilde{s}, \tilde{\theta})$ 
is a local coordinate system on $M_K^2$.   Thus $\diag = \{\rho = \tilde{\rho}, t = \tilde{t}, s = \tilde{s}, \theta = \tilde{\theta}\}$ and 
$$
\fdiag_\rho = \{\rho = \tilde{\rho} = 0, t = \tilde{t}\} \ \  \mbox{and}  \quad \diag_s = \{s = \tilde{s} = 0,  
\rho = \tilde{\rho}, t = \tilde{t}, \theta = \tilde{\theta} \}.
$$
To understand the successive blowups, we use polar coordinates.  Thus first let $R \, \Omega = ( \rho, \tilde{\rho}, t-\tilde{t})$,
so $R \geq 0$ and $\Omega = (\Omega_1, \Omega_2, \Omega_3) \in S^2 $ has $\Omega_1, \Omega_2 \geq 0$. 
The second blowup is around $\diag_s$ which in these coordinates is the set $\{s=\tilde{s} = 0, \Omega = \Omega_*,
\theta = \tilde{\theta})$ where $\Omega_* = (1/\sqrt{2}, 1/\sqrt{2}, 0)$. The portion which lies over $\fdiag_\rho \cap \diag_s$ 
also has $R=0$, so the $S^2$-orthant fiber at $R=0$ is blown up at $\Omega_*$.  Choosing any local coordinate $w$ on $S^2$
with $w=0$ at $\Omega_*$, e.g.\ the projective coordinate $(r-1,\tau)$ below, then the second blowup corresponds to
using polar coordinates $ R' \, \Omega' = (s, \tilde{s}, w, \theta-\tilde{\theta})$ where $R' \geq 0$ and $\Omega'$ lies
in the orthant of $S^4$ where $\Omega'_1, \Omega'_2 \geq 0$. 

Projective coordinates are somewhat more convenient, especially for computing how vector fields transform.
Thus in the first blowup, the projectivization of the $2$-sphere is given by $r = \rho/\tilde{\rho}$ and $\tau = (t-\tilde{t})/\tilde{\rho}$;
the variable $\tilde{\rho}$ serves the place of $r$. These coordinates are valid away from the face $\Omega_2 = 0$,
and together with $\tilde{t}, s, \tilde{s}, \theta, \tilde{\theta}$ give a coordinate system near $\diag \cap \fdiag_\rho$. 
Since $\diag_s = \{r=1, \tau = 0, s = \tilde{s} = 0, \theta = \tilde{\theta}\}$, then for the second blowup we can take 
the projective coordinates on the $4$-sphere orthant, $\sigma = s/\tilde{s}$, $\hat{r} = (r-1)/\tilde{s}$, $\hat{\tau} = \tau/\tilde{s}$, 
$\hat{\theta} = (\theta - \tilde{\theta})/\tilde{s}$.  Observe also that 
$$
\diagiie = \{ r = 1, \tau = 0 \} = \{ \hat{r} = 0, \hat{\tau} = 0, \sigma = 1, \hat{\theta} = 0\},
$$
and this has a clean intersection with the two front faces $\ff_\rho = \{ \tilde{\rho} = 0\}$ and $\ff_s = \{ \tilde{s} = 0\}$. 
The fibers of $\pi_\rho$ are identified with half-spaces $\RR^2_+ = \{ (r,\tau): r \geq 0\}$, while the fibers of
$\pi_s$ are half-spaces $\RR^4_+ = \{ (\sigma, \hat{r}, \hat{\tau}, \hat{\theta}): \sigma \geq 0\}$.  The singular
fiber of $\pi_\rho$ is the union of a copy of $\RR^2_+$ and the blowup $[\RR^4_+; \{(1, 0, 0, 0)\} ]$. 

We have emphasized these singular fibers of $\pi_\rho$ because they are responsible, in some sense, for the new features in the 
analysis over what was done in \cite{MW}. 

\subsubsection{Lifted $\ie$ vector fields, symbols and normal operators}
We now compute the lifts of the generators for the $\ie$ vector fields to $(M_K)^2_{\ie}$ via the left factor, focusing on their 
behavior near $\ff_\rho$ and $\ff_s$.  The simplest formul\ae\ are in projective coordinates. Thus near $\ff_\rho$, in 
the first set of projective coordinates, 
$$
\rho s \del_\rho   =  s r \del_r, \ \rho s \del_t = s r \del_\tau.
$$
The other two vector fields, $s\del_s$, $s\del_\theta$ are already nonvanishing and transform `trivially' in this region.
The key point is that these lifts no longer vanish at the front face $\tilde{\rho} = 0$.  These, together with 
$s\del_s$, $s\del_\theta$, span the normal bundle of $\diagiie$ in this region.   Near $\ff_s$ we compute the lifts all of four generators: 
$$
\rho s \del_\rho =  s r \del_r = \sigma (1 + \tilde{s} \hat{r}) \del_{\hat{r}},\ 
\rho s \del_t = s r \del_\tau = \sigma (1 + \tilde{s} \hat{r}) \del_{\hat{\tau}}, 
s\del_s = \sigma \del_\sigma, s\del_\theta = \sigma \del_{\hat{\theta}}.
$$
These span the normal bundle of $\diagiie$ near $\ff_\rho \cup \ff_s$, which vindicates this blowup as the space which desingularizes this
class of operators.

We have now shown that $\diagiie$ is a $p$-submanifold of $(M_K)^2_\ie$ and that the lifts of the generators of $\calV_\ie$ 
are maximally nondegenerate.  An immediate consequence is that since the operator $\rho s \LKW$ is an `elliptic combination'
of these generating vector fields, its lift is transversely elliptic across $\diagiie$.  This is closely related to the easily checked observation
that the `formal' definition of the $\ie$ symbol ${}^\ie\sigma_1(\rho s \LKW)$ in \eqref{iiesymbol}  agrees with the `ordinary' symbol 
of this lift, i.e., the one obtained by replacing $\del_{\hat{r}}, \del_{\hat{\tau}}, \del_\sigma, \del_{\hat{\theta}}$ by the linear cotangent variables. 

Finally, recall also the formal definition of the normal operators for $\LKW$ at $\ff_\rho$ and $\ff_s$. The former of these is
$$
N_\rho(\LKW) = B_\rho \del_r + B_t \del_\tau + \frac{1}{r}\calJ_0,
$$
where the coefficients are evaluated at $p_0 = (0, t_0) \in K$ and $\calJ_0$ acts on fields over $Z$.  It can be checked immediately
that this is the same as the restriction to $\ff_\rho$ of $\tilde{s} \beta^*(\rho s \LKW)$.  Similarly, $N_s(\LKW)$ 
is the restriction to $\ff_s$ of $\tilde{\rho}\tilde{s} \beta^* (\rho s \LKW)$. 

These observations also exhibit the compatibility between these two normal operators at the corner $\ff_\rho \cap \ff_s$. 

\subsubsection{$\ie$ pseudodifferential operators}
The space of $\ie$ pseudodifferential operators on $M_K$, $\Psi^*_{\ie}(M_K)$, is a direct generalization of 
the space $\Psi^*_0(M)$ of uniformly degenerate pseudodifferential operators given in \cite{Ma, MW}.  
By definition, an operator $A \in \Psi^*_{\ie}(M_K)$ is characterized by the regularity properties of the lift of its Schwartz kernel 
$K_A$ to $(M_K)^2_\ie$, which is required to have a classical singularity across $\diagiie$ and to have 
polyhomogeneous expansions at all boundary faces of $(M_K)^2_\ie$.  

More precisely, $K_A = \beta_* \kappa_A$,  where $\kappa_A$ is a distribution on $(M_K)^2_{\ie}$ which is
a sum of two terms, $\kappa_A^1 + \kappa_A^2$. The first summand, $\kappa_A^1$, is supported near $\diag_\ie$
and vanishes near all boundary faces except the two front faces. It has a classical conormal singularity along $\diagiie$ 
and is the product of specific (dimensionally determined) powers of the boundary defining functions for $\ff_s$ and $\ff_\rho$ 
with a distribution which is smoothly extendible across both of these front faces. The other term, $\kappa_A^2$, 
is smooth in the interior and polyhomogeneous at all boundary faces of $(M_K)^2_\ie$.  The salient features of such an operator
are encoded by several pieces of data: the conormal order of the singularity along the lifted diagonal, which is the order of
the operator, and index sets for the polyhomogeneous expansions at all boundary faces. This leads to somewhat
ornate notation which we for the most part suppress here. 

Before proceeding, we describe the basic example: the identity operator.  In the original 
coordinates, $K_I = \delta(\rho - \tilde{\rho}) \delta(t - \tilde{t}) \delta(s-\tilde{s}) \delta(\theta - \tilde{\theta})$.
Using the two sets of projective coordinates, and recalling the familiar homogeneity of the delta function, we have
$$
\kappa_I = \tilde{\rho}^{-2} \delta( r-1) \delta(\tau) \delta(s-\tilde{s}) \delta(\theta - \tilde{\theta}) = 
\tilde{\rho}^{-2} \tilde{s}^{-4} \delta(\hat{r}) \delta( \hat{\tau}) \delta(\sigma-1) \delta(\hat{\theta}). 
$$
This is a distribution supported on $\diagiie$, so the second type of summand $\kappa_A^2$ does not appear here,
and blows up to order $2$ at $\ff_\rho$ and to order $4$ at $\ff_s$.   This provides a normalization. We demand in general
that the summand $\kappa_A^1$ above equals $R_\rho^{-2} R_s^{-4} \kappa_C$ where $\kappa_C$ has
a classical singularity along $\diagiie$ and is smoothly extendible across $\ff_\rho$ and $\ff_s$, and where
$R_\rho$ and $R_s$ are boundary defining functions for $\ff_\rho$ and $\ff_s$. Any operator
$A$ for which $\kappa_A$ blows up less quickly than this is said to vanish to some order at the front faces. 

There are two key sets of results which make this space of operators a useful tool. The first is a composition law,
which states that the composition of two such operators is again of the same type. This law is accompanied
by some arithmetic which describes the index sets of the composition in terms o the index sets of the 
two factors.  The second describes the boundedness properties of these operators on various function spaces.
In our minimal treatment of this theory, we circumvent all but the most trivial parts of the composition law. 
On the other hand, the boundedness properties are fundamental and we prove them below. 

\subsubsection{The parametrix}
We now turn to the construction of the parametrix $G$.  It is a sum of three terms, $G_0 + G_s + G_\rho$, 
each obtained by solving model problems along $\diagiie$ and the fibers of $\ff_s$ and $\ff_\rho$, respectively. 
In fact, we shall find two such operators, $G_\ell$ and $G_r$, such that $\LKW G_rl = I - R_r$ and $G_\ell \LKW = I - R_ell$,
where $R_r$ and $R_\ell$ are compact.  We shall find such operators where these remainders are smoothing in
the interior and vanish at the front faces; a corollary of the boundedness properties we shall prove below
is that this is enough to ensure that these operators are compact. For simplicty, we focus on the construction of $G_\ell$ 
first, and simply call it $G$. 

The general scheme is that we try to solve the distributional equation $\LKW \kappa_G = \kappa_I$ for
the lifted Schwartz kernel of $G$.  For simplicity we shall simply write these Schwartz kernels as $G$, $I$, etc.
The only singularity of the right hand side is along the lifted diagonal, and the first step is universal in any 
elliptic parametrix construction: we solve away this diagonal singularity using the symbol calculus. 
More precisely, we determine the complete classical expansion of $G_0$ along $\diagiie$ by solving
a sequence of equations using the symbol calculus. This is done on each fiber of the normal bundle of $\diagiie$.
Notably, the construction is uniform in an appropriate sense up to the front faces. We thus choose $G_0$ to 
have support in a small neighborhood of $\diagiie$ not intersecting any boundary face except $\ff_\rho \cup \ff_s$
and to satisfy 
$$
\LKW G_0 = I - R_0,
$$
where $R_0 = R_\rho^{-2} R_s^{-4} \tilde{R}_0$, where $\tilde{R}_0 \in \calC^\infty( (M_K)^2_{\ie})$.  Since the coefficients of $\LKW$ 
blow up like $1/\rho s$,  $G_0$ blows up like $R_\rho^{-1} R_s^{-3}$, i.e., is one order better than the identity operator.

The fact that $R_0$ has the same growth order as $I$ at the front faces means that it is not compact. The next
step in the construction is to choose a correction term $G_\ff = G_\rho + G_s$ such that $\LKW G_\ff = R_0 - R_1$, 
where $R_1$ now vanishes one order faster relative to $R_0$ and hence is compact.  We do this by solving the
model problems
$$
N_s (\LKW) \left. G_s \right|_{\ff_s} = \left. R_0 \right|_{\ff_s}, \qquad \mbox{and} \quad 
N_\rho (\LKW) \left. G_\rho \right|_{\ff_\rho} = \left. R_0 \right|_{\ff_\rho} 
$$ 
along the front faces.  For this to be possible, it is necessary that these normal operators are invertible, or at least surjective.
Since the right hand sides are compactly supported on the fibers of the respective faces, one expects the solutions to be 
polyhomogeneous on these faces. This is carried out exactly as in \cite[Section 5]{Ma}
Compatibility of the normal operators at the corner ensures that these solutions agree there. 
Hence there exists a polyhomogenous Schwartz kernel $G_\ff$ which has the correct leading term at each of these faces.  

Altogether, we have now arranged that
$$
\LKW( G_0 + G_\ff) = I - R_1,
$$
where $R_1$ is smooth in the interior, polyhomogeneous at all boundary faces of $(M_K)^2_{\ie}$ and vanishes to 
order $1$ at $\ff_s \cup \ff_\rho$, hence is a compact operator on certain weighted Sobolev and H\"older spaces
to be described below. 

We have now produced a right parametrix for $\LKW$. An identical construction can be carried out for the adjoint $\LKW^\dag$,
which yields a parametrix $G_\ell^\dag$. Thus $\LKW^\dag G_\ell^\dag = I - R_2^\dag$. Taking transposes gives that
$$
G_\ell \LKW = I - R_2,
$$
i.e., $G_\ell$ is a left parametrix for $\LKW$. 

It is possible to extend this parametrix construction further so that the error term for the right parametrix 
vanishes to infinite order at both $\ff_\rho \cup \ff_s$ and the `left faces', i.e., the lifts of $F_{\rho,\ell} \cup F_{s,\ell}$. 
However, this uses the composition law for $\Psi^{*}_{\ie}$, and we are choosing the simpler route of not proving
this formula.  We can still obtain the some conclusions, but the iteration to obtain higher and higher
regularity is  to a different part of the argument.

\subsubsection{Boundedness}
We now turn to proving boundedness of $\ie$ pseudodifferential operators on weighted H\"older spaces. 
Corresponding results hold for these operators acting between weighted Sobolev spaces, but since these 
are not required here, we omit those proofs. 

We begin with the definition of the $\ie$ H\"older spaces $\calC^{k,\alpha}_\ie(M_K)$. In local coordinates, 
$$
[u]_{\ie; 0, \alpha} = \sup 
\frac{ |u(\rho, t, s, \theta) - u(\rho', t', s' ,\theta')|(\rho + \rho')^\alpha  (s+\s')^\alpha}
{ |\rho - \rho'|^\alpha + |t-t'|^\alpha  + (\rho + \rho')^\alpha(|s-s'|^\alpha + |\theta - \theta'|^\alpha ) } 
$$
and
$$
\calC^{0,\alpha}_{\ie}(M_K) = \{u: ||u||_{\ie; 0, \alpha} = \sup |u| + [u]_{\ie; 0, \alpha} < \infty \}, 
$$
and the associated higher order spaces are
$$
\calC^{k,\alpha}_\ie = \{u:  V_1 \ldots V_j u \in \calC^{0,\alpha}_\ie \ \forall\, V_i \in \calV_{ie},\ j \leq k\}.
$$
More invariantly, these are the natural H\"older spaces associated to the complete metric 
$\hat{g} = \rho^{-2} s^{-2} g$, so for example
$$
||u||_{\ie; 0, \alpha} = \sup |u| + \sup_{p \neq p'} \frac{ |u(p) - u(p')|}{ d_{\hat{g}} (p, p')}.
$$

The important observation is that these norms are scale-invariant in the sense that they remain
unchanged with respect to the two-parameter family of dilations
$$
T_{\lambda, \mu} (\rho, t, s, \theta) =  (\lambda \mu \rho, \lambda \mu t, \mu s, \mu \theta),\ \lambda, \mu > 0.
$$
This is most easily seen using the invariant definition above since these dilations preserve $\hat{g}$ distances (at least
up to bounded factors).

We now state the basic mapping property.
\begin{prop}
Let $A$ be an element of $\Psi^0_{\ie}$ for which its Schwartz kernel has expansions at each of the boundary faces of
the double space, with leading terms at each of these faces
$$
R_{\rho,\ell}^0 \ \mbox{at}\ F_{\rho,\ell},\ \  R_{s,\ell}^\lambda \ \mbox{at}\ F_{s,\ell}, \ \ R_\rho^{-1} \ \mbox{at}\ \ff_\rho,\ \ 
R_s^0\ \mbox{at}\ \ff_s.
$$
and
$$
R_{\rho,r}^0 \ \mbox{at}\ F_{\rho,r},\ \  R_{s,r}^\lambda \ \mbox{at}\ F_{s,r}, \ \ R_\rho^{-1} \ \mbox{at}\ \ff_\rho,\ \ 
R_s^0\ \mbox{at}\ \ff_s.
$$
(We have not stated the optimal exponents here, but simply given ones adequate for the present purposes.)
Then for any $\epsilon < 1$ and $k \geq 0$, 
$$
A: \rho^{-\epsilon} s^{\lambda+ \epsilon} \calC^{k, \gamma}_{\ie}(M_K) \longrightarrow s^\lambda \calC^{k,\gamma}_{\ie}
$$
is bounded.
\end{prop}
\medskip
\noindent{\it Proof:} 
Recall that $A$ decomposes as $A' + A''$ where the Schwartz kernel of $A'$ is supported near $\diagiie$ and
$A''$ has no diagonal singularity. We study these terms separately. 

We first claim that if $u \in \rho^a s^b \calC^{k,\alpha}_\ie$ for any $a, b$, then $A' u \in \rho^{a+1} s^b \calC^{k,\alpha}_\ie$. 
This can be reduced to the local boundedness of uniformly degenerate pseudodifferential operators and then
iteratively to that of standard pseudodifferential operators on ordinary H\"older spaces. First note that 
$\rho^{-a}s^{-b} A' \rho^a s^b$ is an operator of exactly the same type as $A'$, so we may as well assume
that $a=b=0$. Now, decompose the interior of $M_K$ into a countable number of Whitney cubes $Q_\ell$ where 
the $g$-diameter $d_\ell$ of each $Q_\ell$ is one half the distance of $Q\ell$ to $\del_K M_K$ (alternately, this
is a covering by cubes of bounded size with respect to the `partially complete' metric $\rho^{-2} g$). Assuming that
the cover is chosen so that each point is contained in no more than $C$ of these cubes, we have
$$
C_1 \leq ||u||_{\ie; k, \alpha} \leq \sum_\ell || u|_{Q_\ell} ||_{\ie; k, \alpha} \leq C_2 ||u||_{\ie; k, \alpha}. 
$$

Let $D_c$ denote the dilation $(\rho, t, s, \theta) \mapsto (c \rho, c t, s,\theta)$. The idea is based on the
observation that the $\ie$ H\"older norms are invariant with respect to $D_c$ and that the family of
order zero pseudodifferential operators $c^{-1} D_c^* A'$ are uniformly equivalent as $c \to 0$. Writing $D_{d_\ell} = D_\ell$ for 
simplicity, then $D_\ell^{-1} (Q_\ell)$ is a cube of fixed size and of fixed distance from $\del_K M_K$. Therefore, the desired estimate
for $A'$ follows from 
$$
|| u|_{Q_\ell} ||_{\ie; k, \alpha}  = || D_\ell^* (u|_{Q_\ell}) ||_{k,\alpha}
$$
and 
$$
||d_\ell^{-1}D_\ell^* (A' u|_{Q_\ell}) ||_{k,\alpha} \leq C ||D_\ell^* (u |_{Q_\ell}) |_{k,\alpha}.
$$

The estimate for $A''$ proceeds differently. First note that $A'' u$ is $\calC^\infty$ in the interior of
$M_K$; since the vector fields $V \in \calV_{\ie}$ lift to be smooth and tangent to all boundaries
of $M_K$,  $V_1 \cdots V_N A''$ has the same structure as $A''$ for all $N \geq 0$ and $V_i \in \calV_{\ie}$,
hence if we can prove that $|A'' u| \leq C \rho^0 s^\lambda$ then $A'' u \in s^\lambda \calC^{k,\alpha}_\ie$
for all $k \geq 0$.  To obtain this estimate, write the integral expression for $A'' u$ as a pushforward,
$$
A'' u = (\pi_\ell)_* ( \kappa(A'') \pi_r^* u),
$$
where $\pi_\ell, \pi_r$ are the two projections $(M_K^2)_\ie \to M_K$. This can be written out as an integral, but
to do so one would need to write explicit expressions over different coordinate regions in $(M_K^2)_\ie$. In any
case, 
$$
|A'' u| \leq (\pi_\ell)_* ( \kappa(A'') \pi_r^*(\rho^{-\epsilon} s^{\lambda-\epsilon})),
$$
and by the pushforward theorem \cite[Appendix A]{Ma}, this is bounded by $C \rho^0 s^\lambda$.  $\Box$

We also require a closely related result. 
\begin{prop}
If $A$ is an $\ie$ pseudodifferential operator of order $0$ which has the same leading orders at all faces as
in the last Proposition, except at $\ff_\rho$ where it has leading order $R_\rho^{-2}$. Then
$$
A: \calC^{0, \gamma}_{\ie}(M_K) \longrightarrow \log \rho \, \calC^{0,\gamma}_{\ie}
$$
\label{logterms}
\end{prop}
\medskip
\noindent{\it Proof:}    The proof is identical to the previous result except at the very last step. With these hypotheses, 
the fact that $\kappa(A'') \pi_r^* ( \rho^0 s^0) )$ has the same leading coefficient at the abutting faces $F_{\rho,\ell}$ and 
$\ff_\rho$ produces the extra factor of $\log \rho$.  The appearance of this new factor is easy to understand in the following model 
calculation:  consider the pushforward of a smooth compactly supported function $\phi(x,y)$ on the manifold with corners 
$\RR^+_x \times \RR^+_y$ under the map $t = xy$.  Assume that $\phi = 1$ near $(0,0)$ for simplicity. This has a (trivial) expansion 
with leading exponents $0$ at both faces. An elementary change of variables then shows that $\int_{xy = t} \phi = 
-\tfrac12 \log t + \calO(1)$. $\Box$

\subsection{Commutators and higher regularity}
The final part of this analysis uses the parametrix $G$ to deduce higher regularity for solutions of the KW equations
satisfying the generalized Nahm pole boundary conditions.
\begin{theorem} Suppose that $(A,\phi)$ is a solution to ${\KW}(A,\phi) = 0$ on $M$.  Suppose furthermore that
near $K$, $(A,\phi) = (A^K, \phi^K) + (a,\varphi)$, where $(A^K, \phi^K)$ is an approximate solution which is
polyhomogeneous asymptotic to the model knot solution, where $|(a,\varphi)|\leq C \varrho^{-1 + \epsilon} s^{-1 + \varepsilon}$.  
Then $(a,\varphi)$ is polyhomogeneous on $M_K$ with an expansion of the form $a \sim \sum a_0 \varrho^0 + \ldots$, 
$\varphi \sim \varphi_0 \varrho^0 + \ldots$, where all coefficients $a_\gamma$ and $\varphi_\gamma$ in this expansion 
are also polyhomogeneous as $s \to 0$, with $a_\gamma, \varphi_\gamma \sim s^{\alpha} + \ldots$ for some $\alpha > 0$. 
\label{mainregthm}
\end{theorem}
We have not explicitly stated the corresponding regularity theorem for solutions satisfying the Nahm pole condition
at $W \setminus K$; that was the topic of \cite{MW}, and we shall use that result in the following. The proof of this theorem
occupies the rest of this section. 

We begin by writing $\Psi_0 = (A^K, \phi^K)$ and $\Psi = (A, \phi)$, so $\psi = (a,\varphi) = \Psi - \Psi_0$. 
Now rewrite ${\KW}(A,\phi) = 0$ as 
\begin{equation}
\LKW(\psi) =  f + Q(\psi),
\label{Taylor}
\end{equation}
where $f = {\KW}(\Psi_0)$ and $Q(\psi)$ is (precisely) quadratic.    The initial task will be to show that we can choose
the model field $\Psi_0$ so that $f$ vanishes to all orders along $\del M_K$.   We explain this below. 

Granting this, there are several steps to prove that $\psi$ has a polyhomogeneous expansion at $\del M_K$. 
The main part of the argument involves showing that $\psi$ has stable regularity with respect to a a sequence 
of increasingly less degenerate set of vector fields.  
\begin{definition}
Let $\calV$ be any locally finitely generated space of smooth vector fields on a manifold with corners $X$ which is a $\calC^\infty(X)$ 
module and is closed under Lie bracket. We say that a function (or field) $\psi$ has $\calV$-stable regularity, and write
$\psi \in \calA_{\calV}(X)$, if 
$$
V_1 \ldots V_k \psi \in L^\infty(X)\ \mbox{for all}\ V_j \in \calV\ \mbox{and all}\ k.
$$
There are also weighted versions of these spaces, i.e., we can require these arbitrarily many derivatives with respect to elements of
$\calV$ to lie in $\rho_1^{\mu_1} \ldots \rho_N^{\mu_N} L^\infty$, where the $\rho_j$ are boundary defining functions for the
various hypersurface boundaries of $X$ and the $\mu_j$ are fixed weights. 
\end{definition}
One may equally well use any other fixed Banach space in place of $L^\infty$ here. The important point is that {\it all} $\calV$ derivatives
of $\psi$ remain in the same space.   Assuming that the restriction of $\calV$ to the interior of $X$ generates the space of all
smooth vector fields there, then clearly any $\psi \in \calA_{\calV}$ is smooth in the interior. The degeneracies of elements in $\calV$,
however, lead to a well regulated singular structure for $\psi$ at the boundary.

The particular Lie algebras of vector fields we use in the argument below are the following. Each is described in terms of a spanning
set of generators near the boundary, in our standard coordinates.  The first is the space of iterated edge vector fields we have already 
encountered,
$$
\calV_{\ie} = \mathrm{span}\,\{\rho s \del_\rho, \rho s \del_t, s\del_s, s\del_\theta\} ; 
$$
the next is the space of `product $0$' vector fields,
$$
\calV_{p0} = \mathrm{span}\,\{ \rho \del_\rho, \rho\del_t, s\del_s, s\del_t\};
$$
the final one is the most important space of $b$-vector fields
$$
\calV_b(M_K) = \mathrm{span}\,\{\rho \del_\rho, \del_t, s\del_s, \del_\theta\}.
$$
We denote the corresponding spaces of fields by $\calA_{\ie}$, $\calA_{p0}$ and $\calA_b$, respectively; the latter space is simply
called the space of conormal fields.  Observe that $\psi \in \calA_{\ie}$ if and only if $\psi \in \calC^{k,\alpha}_{\ie}$ for every $k$. 
There are spaces of adapted H\"older spaces associated to these other Lie algebras, and the spaces of stable regularity fields
have the similar characterization in terms of these.  There are invariant definitions for the first and third of these spaces of vector fields. 
First, $\calV_b$ is the space of all smooth vector fields on $M_K$ which are tangent to the two boundaries of $M_K$.
On the other hand, $\calV_{\ie}$ is the space of all smooth vector fields which vanish at $s=0$ (the original boundary) and
which lie tangent to the hemisphere fibers along the front face. For $\calV_{p0}$ however, we need to impose some extra
geometric structure near $K$; the vector fields here lie tangent to the fibers of the front face, as for $\calV_{\ie}$, but along
the $s=0$ boundary lie tangent to the $(\rho, t)$ surfaces, which we think of as a local family of surfaces normal to $K$ in $W$. 
Since the result is local near $K$, there is no difficulty in choosing such a family, and we shall abuse notation by thinking
of these vector fields as defined globally. 

The regularity argument proceeds by using the parametrix $G$ constructed above to show first that $\psi \in \calA_{\ie}$,
next that $\psi \in \calA_{p0}$ and finally that $\psi \in \calA_b$.   There is an intermediate step to improve the growth
of $\psi$ from $\calO( (\rho s)^{-1 + \varepsilon}$ to bounded. The passage between these different levels of regularity relies
on the fine mapping properties of $G$, and in particular, in geetting to conormality we require that 
the commutator of $G$ with the $b$ vector fields is better than expected. 

Once we have proved that $\psi$ is conormal, the remaining step to show that it is polyhomogeneous, with exponents at each 
boundary face determined by the indicial roots of $\LKW$ in $\rho$ and $s$, is handled by a formal argument very 
similar to the construction of infinite order approximate solution below. 

We now give details for each of these steps. 

\medskip

\noindent {\bf Construction of the approximate solution $\Psi_0$:}
Recall from Section 4.3 that $(A_0, \phi_0)$ is an approximate solution of the KW equations with Nahm pole boundary condition 
if it equals the model solutions $(A_\varrho, y^{-1} \phi_\varrho)$ and $(A^K, \phi^K)$ at $\wfa$ and $\ff$, respectively, up to lower order 
error terms.  Our immediate goal then is to construct a better approximate solution. We first obtain a formal series solution 
to the equations by iteratively solving for successive terms in this series; taking a Borel sum of this formal series gives a 
field $(A_0, \phi_0)$ which solves the KW equations to infinite order at both boundaries. 

We first carry out the correction steps along $\ff$. This is local at each point $t_0 \in K$, but global on the hemisphere
fiber at that point. Write 
$$
\LKW = B_\rho(0,t_0) + \frac{1}{\rho} \calJ(0, t_0) + E
$$
where we evaluate $B_\rho$ and $\calJ$ along the fiber over $(0, t_0) \in K$. The error term is the collection of all terms
in $\LKW$ for which $\LKW (\rho^\gamma \Psi) = \calO(\rho^\gamma)$, exactly as in \eqref{indface}.  Near this point, write
$\Psi_0 \sim (A^K, \phi^K) + \psi_0 + \rho \psi^1 + \ldots$. Expanding ${\KW}(\Psi_0)$ as in \eqref{Taylor} and inserting
this formal series gives a sequence of equations
$$
(B_\rho(0,t_0) j + \calJ(0,t_0)) \rho^{j-1} \psi_j = f_j \rho^{j-1},
$$
where $f_j$ is the accumulation of all error terms coming from previous steps in the iteration for which the accompanying factors are precisely $\rho^{j-1}$. 
So long as the operator in parentheses on the left is invertible, we can solve this equation for $\psi_j$ and proceed to the next step.
However, if it has nullspace, we can proceed by adding the term $\tilde{\psi}_j \rho^j \log \rho$ to the formal series. In this
process, positive integer powers of $\log \rho$ may accrete as we proceed up the series. but this causes no problems. 
Note also that at such instances, there is a freedom in which solution we choose, i.e., we may add an element of this problematic nullspace.
We may choose this extra term in any reasonable way, but to be systematic we may as well choose it to be $0$.
The solution functions $\psi_j$ are themselves polyhomogeneous at $\del Z = \{s = 0\}$, and some careful bookkeeping
shows that each of them behaves at worst like $\log s$.
 
This procedure is clearly smooth in $t_0$, i.e., the solution functions $\psi_j$ are smooth on $\ff$. We produce in this way 
an infinite formal polyhomogeneous series which solves the equations to infinite order along the front face. 
We choose a Borel sum $\psi'$ for this series. By construction, ${\KW}( \Psi_0 + \psi') = f'$ is a polyhomogeneous 
function on $M_K$ which vanishes to all orders at $\ff = \{ \rho = 0\}$, and has leading term some smooth multiple
of $\log s$ along $\wfa$. 

We can solve away the error term along the original boundary in essentially the same way. For this, at any point $p \in \wfa$, we write
$$
\LKW = B_y(p) \del_y + \frac{1}{y} B_0(p) + E,
$$
and successively solve a sequence of equations $(B_y(p) j + B_0(p)) \psi_j' y^{j-1} = f_j y^{j-1}$. As before, this may require extra factors
of $\log y$. Unlike the previous case, the operators $B_0$ are simply matrices. This is done at each $p \in \wfa$, and the solutions
depend smoothly on this boundary variable.  The solutions do not `spread' on the boundary, and because each $f_j$
vanishes rapidly as $\rho \to 0$, the solutions decay rapidly at $\ff$ as well.    We refer also to \cite{HeMik}, where the
terms in this expansion (in the absence of a knot) are examined closely.

Now choose a Borel sum $\psi''$ for this second asymptotic series.   The complete approximate solution is
$\Psi_0 = \Psi_{\mathrm{model}}  + \psi' + \psi''$, and by construction, ${\KW}(\Psi_0) = f$ vanishes to all orders as $\rho \to 0$
and as $s$ or $y \to 0$. 

\medskip

\noindent {\bf Stable regularity with respect to $\calV_{\ie}$:}
This is an instance of a general argument explained carefully in \cite{ALMP2} and may be applied to solutions of 
degenerate elliptic operators 
associated to any locally finitely generated Lie algebra of vector fields as above. Consider the metric $\hat{g} = \rho^{-2} s^{-2} g$, 
where $g$ is the fixed (incomplete) metric on $M$, lifted to $M_K$. It is straightforward to check that $\hat{g}$ is complete and has 
uniformly bounded geometry, i.e., there is a uniform lower bound for the injectivity radius and the interior of $M_K$ is covered by 
diffeomorphic images of tangent $\epsilon$ balls, for some fixed $\epsilon > 0$, so that the metric coefficients in each of the
associated normal coordinate systems are uniformly bounded in $\calC^\infty$.  The vector fields in $\calV_\ie$ are precisely the 
vector fields on the interior of $M_K$ which have uniformly bounded lengths with respect to $\hat{g}$ and which extend smoothly
to $\del M_K$.  We now invoke Shubin's algebra of uniform pseudodifferential operators, $\Psi^*_{\mathrm{unif}}$ cf.\ \cite{ALMP}.
This exists on any complete manifold $(N,g)$ of bounded geometry. Operators in this class are pseudodifferential operators on $N$ with
Schwartz kernels supported in a uniform neighborhood of the diagonal, $\{ (p,q) \in N^2: \mbox{dist}_{\hat{g}} (p,q) \leq C\}$,
and which are uniformly controlled in all Riemann normal coordinate charts of radius $C$. Compositions of these operators 
and symbol calculus computations carry over immediately from the compact case. Note that the $\ie$ symbol ellipticity
of $\rho s \LKW$ is simply the condition that $\LKW$ is uniformly elliptic in these coordinate charts. The elliptic parametrix construction
yields a left parametrix $\widehat{G} \in \Psi^{-1}_{\mathrm{unif}}$ satisfying
$$
\widehat{G}' \circ (\rho s \LKW) = I - \widehat{R},
$$
where $\widehat{R} \in \Psi^{-\infty}_{\mathrm{unif}}$. Now define $\widehat{G}(\rho, t, s, \theta, \tilde{\rho}, \tilde{t},
\tilde{s}, \tilde{\theta}) = \widehat{G}'(\rho, t , s , \theta, \tilde{\rho}, \tilde{t}, \tilde{s}, \tilde{\theta})  \tilde{\rho} \, 
\tilde{s}$. Applying $\widehat{G}$ to \eqref{Taylor} yields
$$
\psi = \widehat{R} \psi + \widehat{G} f + \widehat{G} Q(\psi),
$$
hence if $V \in \calV_\ie$, then 
$$
V \psi =  V \widehat{R} \psi + V \widehat{G} f + V \widehat{G} Q(\psi).
$$
Now recall that the $\ie$ vector fields are uniform elements of order $1$ in this calculus. Clearly $V \widehat{R} \in 
\Psi^{-\infty}_{\mathrm{unif}}$, so $V \widehat{R} \psi \in (\rho s)^{-1+\varepsilon} \calC^{k,\alpha}_\ie$ for
every $k$.  The second term is even better, since it is smooth and vanishes to infinite order at $\del M_K$.
Hence both of these terms are well controlled. For the third term, note that $V \widehat{G}$ is an operator
of order $0$ with the same expansions at all the boundaries of $(M_K)^2_{\ie}$ as $\widehat{G}$, and furthermore,
$Q(\psi) \in (\rho s)^{-2 + 2\varepsilon} \calC^{2,\alpha}_{\ie}$. Thus $V \widehat{G} Q(\psi) \in (\rho s)^{-1 + 2\epsilon} \calC^{2,\alpha}_{\ie}$.
Altogether, the entire right hand side lives in $(\rho s)^{-1 + \varepsilon} \calC^{2,\alpha}_{\ie}$, and this is true for any
$V \in \calV_{\ie}$, so $\psi \in (\rho s)^{-1 + \varepsilon} \calC^{3,\alpha}_\ie$. Iterating this argument shows that
$\psi \in (\rho s )^{-1 + \varepsilon} \calC^{k,\alpha}_{\ie}$ for every $k$, so $\psi$ has stable regularity with
respect to $\calV_{\ie}$. 

\medskip

\noindent{\bf From $(\rho s)^{-1 + \varepsilon} \calA_{\ie}$ to $\calA_{\ie}$}
The previous argument did not allow us to obtain any improvement in the growth rate of $\psi$. For this we must use
the refined parametrix $G$.  The key advantage is that both $G$ and $R$ vanish to order one at $\ff_\rho$ and $\ff_s$,
so for example, $G$ maps $(\rho s)^\lambda\calC^{k,\alpha}_\ie$ to $(\rho s)^{\lambda + 1} \calC^{k +1,\alpha}_\ie$ if $\lambda < -1$,
and to $\calC^{k+1,\alpha}_\ie$ if $\lambda \geq -1$.  One cannot gain a factor of $\rho s$ if $\lambda > -1$ because 
$G$ and $R$ are only bounded at the left face. 

We apply this to 
\begin{equation}
\psi = R \psi + G f + G Q(\psi),
\label{theday}
\end{equation}
recalling that $\psi \in (\rho s)^{-1 + \varepsilon} \calC^{k,\alpha}_\ie$ for every $k$. The first term on the right lies in $\calC^{k,\alpha}_\ie$,
the second term is smooth and vanishes to all orders, and since $Q(\psi) \in (\rho s)^{-2 + 2\varepsilon} \calC^{k,\alpha}_\ie$,
the third term lies in $(\rho s)^{-1 + 2\varepsilon} \calC^{k,\alpha}_\ie$.  Thus we have gained a factor of $(\rho s)^{\varepsilon}$.

Iterating this argument a finite number of times shows that $\psi \in \calC^{k,\alpha}_\ie$ for every $k$. 

\medskip

\noindent {\bf From $\calA_{\ie}$ to $\calA_{p0}$:} 
To improve $\ie$ stable regularity to stable regularity with respect to $\calV_{p0}$,  we must control derivatives
with respect to $\rho \del_\rho$ and $\rho \del_t$, rather than just $\rho s \del_\rho$ and $\rho s \del_t $. 
Apply either one of these derivatives to \eqref{theday}; for example, 
$$
\rho \del_\rho \psi = \rho \del_\rho R \psi + \rho\del_\rho Gf + \rho\del_\rho G Q(\psi).
$$
As noted above, both $G$ and $R$ vanish to order $1$ at $\ff_s$, so $\rho\del_\rho R$ and $\rho\del_\rho G$ are of order $0$ at $\ff_s$
while still of order $1$ at $\ff_\rho$. These compositions are thus bounded on $\calC^{k,\alpha}_\ie$ for all $k$, which implies that each 
of the terms on the right lie in these spaces. The argument for $\rho \del_t$ is similar. 

This proves that $\rho \del_\rho \psi,  \rho \del_t \psi \in \calC^{k,\alpha}_\ie$ for every $k$. 
An induction leads to the conclusion that
$$
(\rho \del_\rho)^i (\rho \del_t)^j (s\del_s)^m (s\del_\theta)^\ell (a, \varphi) \in \calC^{k,\alpha}_\ie
$$
for all $i, j, m, \ell \geq 0$. 

\medskip

\noindent {\bf From $\calA_{p0}$ to conormality:} 
This step relies on the 
\begin{lemma}
If $A \in \Psi^{r, \calE}_\ie$, then each of the commutators 
$$
[\rho \del_\rho, A],\ [\rho \del_t, A],\ [\del_t, A],\ [\del_\theta, A]
$$
are also $\ie$ pseudodifferential operators with the same pseudodifferential order and the same
index sets at all faces of $(M_K)^2_\ie$. 
\end{lemma}
\noindent{\it Proof:} The analogous result in the simple edge case already appears in \cite{Ma}, so we explain 
how to adapt that proof to this slightly more general setting.  Each of the individual operators
$\rho\del_\rho A$, $A \rho \del_\rho$, etc.\, is precisely one order more singular at $\ff_\rho$ and/or $\ff_s$, 
but the leading order singularities cancel when taking the differences in the commutators.   Near $\ff_s$ we
can immediately apply the same proof as in \cite{Ma}, even uniformly up to $\ff_\rho \cap \ff_s$ since
all of these vector fields are tangent to this intersection. 
$\Box$

\medskip

We may now prove the conormality of $\psi$. The new consideration is to control derivatives with respect to $\del_t$ and $\del_\theta$. 
This relies on two facts. First, although $\del_t G$ and $\del_\theta G$ are still of order $0$ at $\ff_s \cup \ff_\rho$, they no 
longer vanish there, and the same is true for $R$. Second, the commutators $[\del_t, G]$, $[\del_\theta, G]$, 
$[\del_t, R]$, $[\del_\theta, R]$ still vanish to order $1$ at $\ff_s \cup \ff_\rho$. 

Now write
$$
\del_t \psi = \del_t R \psi + \del_t G f + \del_t G Q(\psi).
$$
By Lemma \eqref{logterms}, the first and third terms lie in $\log \rho \calC^{k,\alpha}_\ie$ and the second is
rapidly vanishing.  Now write $\del_t R \psi = R \del_t \psi + [\del_t , R ] \psi$.  The operator $R$ maps
$\log \rho \calC^{k,\alpha}_{\p0}$ to $\calC^{k,\alpha}_{\p0}$, and by the lemma above, the third term also 
lies in the same space.  On the other hand,
$$
\del_t G Q(\psi) = G \del_t Q(\psi) + [\del_t, G] Q(\psi) = G  B(\psi , \del_t \psi) + [\del_t, G] Q(\psi).
$$
Here $B$ is a bilinear form with polyhomogeneous coefficients in its two arguments. Certainly 
$B(\psi, \del_t\psi) \in \log \rho \calC^{k,\alpha}_{\p0}$, so applying $G$ to it yields something
in $\calC^{k,\alpha}_{\p0}$. The final term is obviously in $\calC^{k,\alpha}_{\p0}$ as well. 

These same arguments work for $\del_\theta \psi$, and as before we can iterate this
argument to control derivatives of all orders.   Examining the argument above and recalling
the leading powers of $G$ at each of its faces, this argument actually proves that
$\psi \in s^\varepsilon \calC^{k,\alpha}_b$ for every $k$. 

\medskip

\noindent {\bf Polyhomogeneity:}
We now explain how to pass from conormality to the existence of expansions at the two boundary faces of $M_K$. 
The main point is that all tangential derivatives may be treated as lower order. 

We first regard the equation ${\KW}(\Psi_0 + \psi) = 0$ as a nonlinear ODE in $y$, with the dependence on the boundary
variable $z$ purely parametric.  It is then a classical argument that an exact solution of such an ODE has
an expansion in $y$, and it follows from their definition that the exponents in this expansion are regulated
by the indicial roots at this face; this part of the argument is identical to the one in \cite{MW}.  We follow
this method to show that $\psi$ has an expansion in $y$, which near $\ff \cap \wfa$ we regard
as an expansion in the variable $s$.  The dependence of the solution as a function of $t, z$ (or near
the corner, as a function of $\rho, t, \theta$), and hence each coefficient in this asymptotic expansion
is conormal in these variables.

On the other hand, as $\rho \to 0$ we consider ${\KW}(\Psi_0 + \psi) = 0$ as a nonlinear conic elliptic equation on 
$S^2_+ \times \RR^+$ with all data depending smoothly on the parameter $t$.  To say that this is conic
problem means that it takes the form $\del_\rho \psi + \rho^{-1} \calJ \psi  = f + Q(\psi)  + E\psi$,
where $Q(\psi)$ incorporates the quadratic terms and $E\psi$ contains the linear terms involving derivatives that
we are regarding as parametric, i.e., the $t$ derivatives.  This too may be solved iteratively to show
that $\psi$ has a complete expansion as a function of $\rho$ with exponents determined by the indicial roots here,
which are determined by the eigenvalues of $\calJ$. The coefficients are conormal as $s \to 0$ and depend smoothly on $t$. 

To show that these expansions at each of the two faces fit together as a polyhomogeneous expansion on the
manifold with corners $M_K$, which means in particular that there is a product-type expansion at the corner, 
we invoke a result from \cite[Appendix]{Ma}: 
\begin{lemma}
Suppose $\psi$ is a conormal distribution on a manifold with corners $X$ such that at every boundary hypersurface $H$ 
of $X$, $\psi \sim \sum R_H^{\gamma} \psi_{H,\gamma}$, where $R_H$ is the boundary defining function for $H$. Suppose also
that these coefficients are themselves conormal as distributions on $H$ (qua manifold with boundary). 
Then $\psi$ is polyhomogeneous on $X$. 
\end{lemma}

This now concludes the proof that the gauged solutions to ${\KW}(A,\phi)$ which satisfy the generalized Nahm pole
boundary condition with a knot singularity are polyhomogeneous on $M_K$. This is the key technical statement
needed for various arguments earlier in this paper. 

\medskip

\noindent{\it Acknowledgements:} Research of RM supported in part by NSF grant DMS-1608223.   Research of EW supported in part
by NSF Grant PHY-1606531.
The authors are grateful to Siqi He, who read the manuscript carefully and gave some useful remarks.

\bibliographystyle{unsrt}

\end{document}